\UseRawInputEncoding 

\documentclass[10pt]{article}
\usepackage[letterpaper, left=1in, top=1in, right=1in, bottom=1in, verbose, ignoremp]{geometry}

\usepackage{url}\RequirePackage[colorlinks,citecolor=blue, linkcolor=blue,urlcolor = blue]{hyperref}
\usepackage{latexsym,amssymb,amsmath,amsfonts,graphicx,color,fancyvrb,amsthm,enumerate,subcaption,mathrsfs}
\usepackage[longnamesfirst,authoryear,round]{natbib}
\usepackage[dvipsnames]{xcolor}
\usepackage{xy}\xyoption{all} \xyoption{poly} \xyoption{knot}
\usepackage{float}
\usepackage{bm}
\usepackage{bbm}
\usepackage{makecell}
\usepackage{multicol,multirow}
\usepackage{array}
\usepackage{relsize}
\usepackage{chngcntr}
\usepackage{etoolbox}
\usepackage{caption}
\usepackage{tikz}
\usetikzlibrary{decorations.pathreplacing,calc}
\usetikzlibrary{shapes,backgrounds}
\usetikzlibrary{patterns}
\usetikzlibrary{cd}
\usepackage{amsopn}
\usepackage{calc}
\usepackage{algorithm}
\usepackage[noend]{algpseudocode}
\usepackage{hyperref}
\usepackage{mathtools}

\usepackage{tikz-cd} 
\usepackage{amscd} 
\usepackage{comment}
\usepackage[capitalize,nameinlink,compress]{cleveref}
\crefname{figure}{Figure}{Figures} 
\crefname{equation}{}{} 
\crefname{assumption}{Assumption}{Assumptions}
\crefname{subsection}{Subsection}{Subsections}

\newcounter{cdrow}

\newtheorem{theorem}{Theorem}[]
\newtheorem*{theorem*}{Theorem}
\newtheorem{corollary}[theorem]{Corollary}
\newtheorem{lemma}[theorem]{Lemma}
\newtheorem{proposition}[theorem]{Proposition}
\newtheorem{assumption}[theorem]{Assumption}
\newtheorem{claim}[theorem]{Claim}
\newtheorem*{claim*}{Claim}

\theoremstyle{definition}
\newtheorem{definition}[theorem]{Definition}
\newtheorem*{definition*}{Definition}

\theoremstyle{remark}
\newtheorem{remark}[theorem]{Remark}

\newtheorem*{example*}{Example}

\def\log{{\rm log}}
\def\tr{{\rm tr}}

\makeatletter
\newcommand*{\op}{%
	\DOTSB
	\mathop{\vphantom{\bigoplus}\mathpalette\matt@op\relax}%
	\slimits@
}
\newcommand\matt@op[2]{%
	\vcenter{\m@th\hbox{\resizebox{\widthof{$#1\bigoplus$}}{!}{$\boxplus$}}}%
}
\makeatother

\newcommand{\R}{\mathbb{R}}
\newcommand{\N}{\mathbb{N}}

\newcommand{\argmin}{\text{argmin}}
\newcommand{\TPT}[1]{ \R^{#1}/\R \mathbf{1}}
\newcommand{\dtr}{d_{\mathrm{tr}}}

\newcommand{\sgn}{\text{\sgn}}

\renewcommand{\vec}[1]{\mathbf{#1}}

\DeclareMathOperator{\tconv}{tconv}
\DeclareMathOperator{\conv}{conv}

\makeatletter
\def\@biblabel#1{}
\makeatother

\makeatletter
\patchcmd{\NAT@citex}
{\@citea\NAT@hyper@{%
		\NAT@nmfmt{\NAT@nm}%
		\hyper@natlinkbreak{\NAT@aysep\NAT@spacechar}{\@citeb\@extra@b@citeb}%
		\NAT@date}}
{\@citea\NAT@nmfmt{\NAT@nm}%
	\NAT@aysep\NAT@spacechar\NAT@hyper@{\NAT@date}}{}{}

\patchcmd{\NAT@citex}
{\@citea\NAT@hyper@{%
		\NAT@nmfmt{\NAT@nm}%
		\hyper@natlinkbreak{\NAT@spacechar\NAT@@open\if*#1*\else#1\NAT@spacechar\fi}%
		{\@citeb\@extra@b@citeb}%
		\NAT@date}}
{\@citea\NAT@nmfmt{\NAT@nm}%
	\NAT@spacechar\NAT@@open\if*#1*\else#1\NAT@spacechar\fi\NAT@hyper@{\NAT@date}}
{}{}

\makeatother

\begin{document}
	
	\def\spacingset#1{\renewcommand{\baselinestretch}%
		{#1}\small\normalsize} \spacingset{1}

	\begin{flushleft}
		{\Large{\textbf{Tropical Gradient Descent}}}
		\newline
		\\
		Roan Talbut$^{1,\dagger}$ and Anthea Monod$^{1}$
		\\
		\bigskip
		\bf{1} Department of Mathematics, Imperial College London, UK
		\\
		\bigskip
		$\dagger$ Corresponding e-mail: r.talbut21@imperial.ac.uk
	\end{flushleft}
	
	
	\section*{Abstract}
	
	We propose a gradient descent method for solving optimization problems arising in settings of tropical geometry---a variant of algebraic geometry that has {attracted growing interest} in applications such as computational biology, economics, and computer science.  Our approach takes advantage of the polyhedral and combinatorial structures arising in tropical geometry to propose a versatile {method} for approximating local minima in tropical statistical optimization problems---a rapidly growing body of work in recent years.  Theoretical results establish global solvability for 1-sample problems and a convergence rate matching classical gradient descent. Numerical experiments demonstrate the method's superior performance {compared to classical gradient descent} for tropical optimization problems which exhibit tropical convexity but not classical convexity. We also demonstrate the seamless integration of tropical descent into advanced optimization methods, such as Adam, offering improved overall {accuracy}.
	
	\paragraph{Keywords:} tropical quasi-convexity; tropical projective torus; gradient descent; phylogenetics.
	
\section{Introduction} \label{sec:intro}

Where algebraic geometry studies geometric properties of solution sets of systems of multivariate polynomials, \emph{tropical geometry} restricts to the case of polynomials defined by linearizing operations, where the ``sum'' of two elements is their maximum and the ``product'' of two elements is their sum.  Evaluating polynomials with these operations results in piecewise linear functions.  Tropical geometry is a relatively young field of pure mathematics established in the 1990s (with roots dating back further) and has recently become an area of active interest in computational and applied mathematics for its relevance to the statistical analyses of phylogenetic trees \citep{monod2018StatPerspective}, the problem of dynamic programming in computer science \citep{maclagan2015introduction}, and in mechanism design for two-player games in game theory \citep{LIN2019133}.  Following the identification of the tropical Grassmannian and the space of phylogenetic trees \citep{speyer2004grassmannian} in particular, there has been {a} surge of research \citep{barnhill2023tml,lin2018FermatWeber,comuaneci2023transportation,monod2018StatPerspective} which necessitates the development of new computational techniques for {optimization} in tropical geometric settings.  

{Training machine learning models and many computational tasks} in statistics entail solving optimization problems; this paper is motivated by statistical optimization problems in tropical geometric settings. {Early research in tropical statistics focused on} optimization problems such as the identification of Fermat--Weber points which, due to the piecewise linear nature of tropical geometry, can be re-framed as a linear program \citep{comuaneci2023transportation,lin2018FermatWeber}. Unfortunately, the complexity of these linear programs often scales {poorly} with the size of our dataset, so gradient {methods are} more commonly implemented in practice \citep{barnhill2024polytropes,barnhill2023tml}. More recently---following the identification of ReLU neural networks and tropical rational functions \citep{zhang2018tropical}---more sophisticated neural network optimization problems are being studied in the tropical setting \citep{pasque2024tropDBs,yoshida2023tropNNs}, which necessitates the use of gradient methods for optimization in the tropical setting. In this paper, we address the natural question---given the need for gradient methods in tropical data science, how can we tailor classical gradient methods to the tropical setting? We provide a foundational tropical gradient descent framework to solve the optimization problems which arise in the computation of a wide range of tropical statistics, an area of active interest in applied tropical geometry.

The remainder of this paper is organized as follows. The following section contains the necessary foundations of tropical geometry, an introduction to tropical location problems \citep{comuaneci2024location}, and {an overview of the} statistical optimization tasks studied in this paper. \Cref{sec:contribution} contains the theoretical contributions of this paper; we formulate steepest descent with respect to the tropical norm and present the theoretical guarantees for tropical descent, most notably its convergence for a wide class of tropically quasi-convex functions and a convergence rate of {$O(1/\sqrt{m})$}. In \Cref{sec:experiments}, we perform a comprehensive numerical study of our proposed tropical descent method and its performance in the statistical optimization problems of interest. We conclude with a discussion of possible future work on tropical gradient methods in \Cref{sec:conclusion}.

\section{Background and Preliminaries}

In this section, we {present a summary of} tropical geometry and the statistical optimization problems which we study in the tropical setting. 

\subsection{Tropical Geometry}

{We begin by reviewing} the algebraic and geometric structure of the tropical projective torus, the state space for phylogenetic data and the domain of our tropical statistical optimization problems.

\begin{definition}[Tropical Algebra]
	The \emph{tropical algebra} is the semiring $\overline\R = \R \cup \{- \infty\}$ with the addition and multiplication operators---tropical addition and tropical multiplication, respectively---given by
	\begin{align*}
		a \boxplus b = \max \{a,b\}, \quad a \odot b = a+b. 
	\end{align*}
	The additive identity is $-\infty$ and the multiplicative identity is $0$. Tropical subtraction is not defined; tropical division is given by classical subtraction.
\end{definition}

The tropical algebra defined above is sometimes referred to as the \emph{max-plus algebra}; we can define the min-plus algebra analogously, taking $\min \{ a,b \} =: a \oplus b$ as the additive operation on $\R \cup \{ \infty \}$. These are algebraically equivalent {under} negation. Unless specified, we use the max-plus convention as this is better suited for phylogenetic statistics on the tropical projective torus \citep{speyer2004grassmannian}.

\begin{definition}[Tropical Projective Torus]
	The $N-1$-dimensional \emph{tropical projective torus} is a quotient space constructed by endowing {$\R^{N}$} with the equivalence relation 
	\begin{equation}
		\label{eq:tpt_equiv}
		\vec x \sim \vec y \Leftrightarrow \exists \, a { \, \in \R} : \: \vec x = a \odot \vec y;
	\end{equation}
	it is denoted by $\TPT{ N}$.  
	The generalized Hilbert projective metric, also referred to as the \emph{tropical metric}, is given by
	\[
	d_{\tr}(\vec x,\vec y) = \max_i {(}x_i-y_i{)} - \min_i {(}x_i - y_i{)} = \max_{i,j} {(}x_i - y_i - x_j + y_j{)}.
	\]
	This metric is induced by the tropical norm, which is given by
	\[
	\| \vec x\|_{\tr} = \max_i x_i - \min_i x_i.
	\]
\end{definition}

\begin{remark}
	While a single point in the tropical projective torus is given by some ray $\vec x + \R\mathbf{1}$, we note that there is a unique representative of this equivalence class whose coordinates sum to zero. In taking such representatives, we can identify the tropical projective torus $\TPT{ N}$ and the hyperplane $\mathcal{H} = \{\vec x \in \R^{ N} : \sum x_i = 0  \}$. Throughout this paper, we use this homeomorphism to {visualize} the tropical projective torus (e.g., \Cref{fig:TropicalBalls,fig:tropicalgeometry}).
\end{remark}

The tropical projective torus is the ambient space containing the tropical Grassmannian, which is equivalent to the space of phylogenetic trees \citep{speyer2004grassmannian}. This has sparked the study of tropical data science, formalizing statistical techniques which respect the tropical geometry of the state space \citep{monod2018StatPerspective,yoshida2020science}.

The generalized Hilbert projective metric defined above is widely accepted and preferred for the development of geometric statistical tools on the tropical projective torus \citep{yoshida2020science}. However, recent work has shown theoretical advantages to the use of an asymmetric metric \citep{comuaneci2023transportation}.

\begin{definition}[Tropical Asymmetric Distance]
	The \emph{tropical asymmetric distances} on $\TPT{ N}$ are given by
	\begin{align*}
		d_{\triangle_{\min}}(\vec a, \vec b) \coloneqq \sum_i (b_i-a_i) - {N}\min_j(b_j-a_j),\\
		d_{\triangle_{\max}}(\vec a, \vec b) \coloneqq {N}\max_j(b_j-a_j) - \sum_i (b_i-a_i).
	\end{align*}
\end{definition}

\begin{remark}
	The asymmetric distances act as an $\ell_1$ norm on the tropical projective torus, while the symmetric metric $\dtr$ acts as an $\ell_{\infty}$ norm. We note that ${N} d_{\tr} = d_{\triangle_{\min}}+d_{\triangle_{\max}}$. \Cref{fig:TropicalBalls} shows the balls for each tropical metric of interest.
\end{remark}

\begin{figure}
	\centering
	\begin{tikzpicture}[scale = 0.8]
		\coordinate (A) at (0,0);
		\fill[black] (A) circle (2pt);
		\node (A) at (0,0) [above right] {\tiny$\mathbf{0}$};
		\draw [fill=gray,opacity=1/4] (0.87, -0.5) -- (0.87,0.5) -- (0,1) -- (-0.87,0.5) -- (-0.87,-0.5) -- (0,-1) -- (0.87, -0.5);
		\draw[-stealth] (0,0) -- (1.16,0);
		\draw[-stealth] (0,0) -- (-0.58,1);
		\draw[-stealth] (0,0) -- (-0.58,-1);
		
		\coordinate (B) at (5,0);
		\fill[black] (B) circle (2pt);
		\node (B) at (5,0) [above right] {\tiny$\mathbf{0}$};
		\draw [fill=gray,opacity=1/4] (5.87, 0) -- (4.565,0.75) -- (4.565,-0.75) -- (5.87, 0);
		\draw[-stealth] (5,0) -- (6.16,0);
		\draw[-stealth] (5,0) -- (4.42,1);
		\draw[-stealth] (5,0) -- (4.42,-1);
		
		\coordinate (C) at (10,0);
		\fill[black] (C) circle (2pt);
		\node (C) at (10,0) [above right] {\tiny$\mathbf{0}$};
		\draw [fill=gray,opacity=1/4] (9.13, 0) -- (10.435,-0.75) -- (10.435,0.75) -- (9.13, 0);
		\draw[-stealth] (10,0) -- (11.16,0);
		\draw[-stealth] (10,0) -- (9.42,1);
		\draw[-stealth] (10,0) -- (9.42,-1);
	\end{tikzpicture}
	\caption{The $\dtr$, $d_{\triangle_{\min}}$ and $d_{\triangle_{\max}}$ tropical unit balls in $\TPT{3} \cong \mathcal{H} = \{ \sum x_i = 0  \}$. {The solid lines show coordinate directions.}}
	\label{fig:TropicalBalls}
\end{figure}

We now define the primary geometric objects of interest in the tropical projective torus; hyperplanes, lines, and convex hulls. We use the $\max$-convention for these definitions, though we note each has a $\min$-convention equivalent. 

\begin{definition}[Tropical Hyperplane]
	The \emph{tropical hyperplane} $\mathcal{H}_{\vec a}$ defined by the linear form $a_1 \odot x_1 \boxplus \cdots \boxplus {a_N \odot x_N}$ is given by
	\[
	\mathcal{H}_{\vec a} = \{ \vec x { \,\in \TPT{N}} : \exists \, i \neq j\text{ s.t. }a_i+x_i=a_j+x_j = \max_k {(}a_k + x_k{)}  \}.
	\]
\end{definition}

\begin{definition}[Tropical Line Segment]
	For any two points $\vec a, \vec b \in \TPT{ N}$, the \emph{tropical line segment} between $\vec a$ and $\vec b$ is the set
	\[
	\gamma_{\vec a \vec b} = \{ \alpha \odot \vec a \boxplus \beta \odot \vec b \mid \alpha, \beta \in \R\}
	\]
	with tropical addition taken coordinate-wise.
\end{definition}

Tropical line segments define a unique geodesic path between any two points, but general geodesics are not uniquely defined on the tropical projective torus.

\begin{figure}[!ht]
	\begin{subfigure}{0.3\textwidth}
		\hfill
		\centering
		\begin{tikzpicture}[scale = 0.8]
			\node (a) at (0,0) [above right] {\tiny$-\vec a$};
			\coordinate (A) at (0,0);
			\draw [dashed, ->] (A) -- ($(A)+1.75*(-0.577,-1)$);
			\draw [dashed, ->] (A) -- ($(A)+4*(0.5,0)$);
			\draw [dashed, ->] (A) -- ($(A)+1.75*(-0.577,1)$);
			\draw (A) -- ($(A)-1.75*(-0.577,-1)$);
			\draw (A) -- ($(A)-4*(0.5,0)$);
			\draw (A) -- ($(A)-1.75*(-0.577,1)$);
			\fill[black] (A) circle (2pt);
		\end{tikzpicture}
		\label{fig:tropicalhyperplane}
	\end{subfigure}
	\hfill
	\begin{subfigure}{0.3\textwidth}
		\centering
		\begin{tikzpicture}[scale = 0.8]
			\node (B) at (1,2) [above right] {\tiny$\vec b$};
			\node (A) at (-1,0) [left] {\tiny$\vec a$};
			\coordinate (B) at (1,2);
			\coordinate (A) at (-1,0);
			\draw [dashed, ->] (B) -- ($(B)+3*(-0.577,-1)$);
			\draw [dashed, ->] (A) -- ($(A)+1*(-0.577,-1)$);
			\draw [dashed, ->] (B) -- ($(B)+(0.5,0)$);
			\draw [dashed, ->] (A) -- ($(A)+(2.5,0)$);
			\draw [dashed, ->] (B) -- ($(B)+0.5*(-0.577,1)$);
			\draw [dashed, ->] (A) -- ($(A)+1*(-0.577,1)$);
			\draw (A) -- ($(A)+(0.923,0)$);
			\draw (B) -- ($(B)+2*(-0.577,-1)$);
			\fill[black] (B) circle (2pt);
			\fill[black] (A) circle (2pt);
		\end{tikzpicture}
		\label{fig:tropicalline}
	\end{subfigure}
	\hfill
	\begin{subfigure}{0.3\textwidth}
		\centering
		\begin{tikzpicture}[scale = 0.8]
			\node (A) at (1,2) [above right] {\tiny$\vec a$};
			\node (B) at (1.25,-0.5) [above right] {\tiny$\vec b$};
			\node (C) at (-1,0) [left] {\tiny$\vec c$};
			\coordinate (A) at (1,2);
			\coordinate (B) at (1.25,-0.5);
			\coordinate (C) at (-1,0);
			\draw [dashed, ->] (A) -- ($(A)+3*(-0.577,-1)$);
			\draw [dashed, ->] (B) -- ($(B)+0.5*(-0.577,-1)$);
			\draw [dashed, ->] (C) -- ($(C)+1*(-0.577,-1)$);
			\draw [dashed, ->] (A) -- ($(A)+(1.25,0)$);
			\draw [dashed, ->] (B) -- ($(B)+(1,0)$);
			\draw [dashed, ->] (C) -- ($(C)+(3.25,0)$);
			\draw [dashed, ->] (A) -- ($(A)+0.5*(-0.577,1)$);
			\draw [dashed, ->] (B) -- ($(B)+3*(-0.577,1)$);
			\draw [dashed, ->] (C) -- ($(C)+1.5*(-0.577,1)$);
			\draw (C) -- ($(C)+(1.923,0)$);
			\draw (B) -- ($(B)+2*(-0.577,1)$);
			\draw (A) -- ($(A)+2*(-0.577,-1)$);
			\draw [fill=gray,opacity=1/4] (0.923,0) -- (-0.154,0) -- (0.3845,0.95);
			\fill[black] (A) circle (2pt);
			\fill[black] (B) circle (2pt);
			\fill[black] (C) circle (2pt);
		\end{tikzpicture}
		\label{fig:tropicalconvexhull}
	\end{subfigure}
	\hfill
	\caption{A tropical hyperplane, tropical line segment and tropical convex hull in $\TPT{3} \cong \mathcal{H} = \{\vec x:  \sum x_i = 0  \}$. The dashed lines show coordinate directions.}\label{fig:tropicalgeometry}
\end{figure}

\begin{definition}[Tropical Convex Hull \citep{develin2004tropical}]
	A set $S \subset \TPT{N}$ is \emph{max-tropically convex} if it contains $\alpha \odot \vec x \boxplus \beta \odot \vec y$ for all $\vec x, \vec y \in S$ and $\alpha, \beta \in \R$. For a finite subset $X = \{\vec x_1,\ldots, \vec x_{{K}}\} \subset \TPT{{N}}$, the \emph{max-plus tropical convex hull} of $X$ is the set of all max-tropical linear combinations of points in $X$,
	\[
	\tconv_{\max}(X) \coloneqq \{ \alpha_1 \odot \vec x_1 \boxplus \cdots \boxplus \alpha_{ K} \odot \vec x_{ K} \mid \alpha_1, \ldots, \alpha_{ K} \in \R\}.
	\]
	Similarly, the \emph{min-tropical convex hull} of $X$, $\tconv_{\min}(X)$, is the set of all min-tropical linear combinations of points in $X$.
\end{definition}

A study of tropically convex functions requires an understanding of both max-tropical and min-tropical convexity, as we will see in our following consideration of tropical location problems.

\subsection{Tropical Location Problems}

Tropical location problems are a new avenue of research, motivated by the optimization of statistical loss functions; they involve the minimization over some dataset of some loss function which heuristically increases with distance from some kernel. Here we outline the relevant definitions of tropical location problems and their theoretical behavior as presented by \cite{comuaneci2024location}. We introduce the problem of tropical linear regression \citep{akian2023tropical} as motivation.

\subsubsection{Motivating Example}

The problem of tropical linear regression looks to find the best-fit tropical hyperplane which minimizes the maximal distance to a set of data points. While this is an inherently statistical optimization problem, it has recently been proven to be polynomial-time equivalent to solving mean payoff games \citep{akian2023tropical}.

A tropical hyperplane is uniquely defined by a cone point $\vec t \in \TPT{N}$, allowing us to formulate the best-fit tropical hyperplane as an optimization problem over ${\vec t \in \, }\TPT{N}$.

\begin{definition}[Tropical Linear Regression \citep{akian2023tropical}]
	\label{def:lin_reg}
	Let $X = \{\vec x_1, \dots, \vec x_K \} \subset \TPT{N}$ be our dataset. The \emph{tropical linear regression problem} finds the vertex $\vec t \in \TPT{N}$ of the best-fit hyperplane:
	\begin{align*}
		\min_{\vec t \in \TPT{N}} f(\vec t) = \max_{k \leq K} \dtr(\vec x_k, \mathcal{H}_{\vec t}) =  \max_{k \leq K} { \left( \textstyle \max_i (\vec x_k-\vec t)_i - \text{2ndmax}_i (\vec x_k - \vec t)_i  \right)}
	\end{align*}
	{where $\text{2ndmax}_i y_i$ denotes some $y_i$ such that $ \exists \, j \neq i$ such that $\forall \,k \neq i,j$: $y_j \geq y_i \geq y_k$.}
\end{definition}

{This expression for the tropical distance between $\vec x$ and $\mathcal H_{\vec t}$ was first stated by \cite{gartner2008tropical}, and has been frequently used in the ongoing study of tropical support vector machines \citep{yoshida2023tropical,tang2024tropical}.}

\begin{figure}[!ht]
	\centering
	\includegraphics[width = 0.5\textwidth]{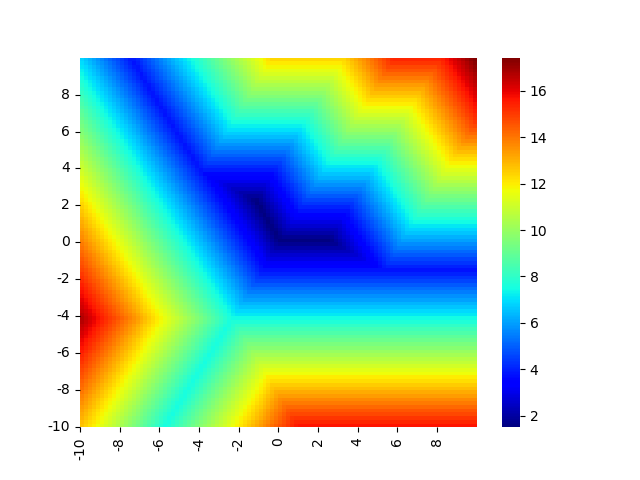}
	\caption{A heat map of a 3-sample linear regression objective on $\TPT{3} \cong \mathcal{H} = \{ \vec x: \sum x_i = 0 \}$.}
	\label{fig:LR_Visualisation}
\end{figure}

In \Cref{fig:LR_Visualisation}, we see a heatmap of this objective function over the tropical projective plane for a sample of size 3. We see immediately that the objective is not convex. The linear regression objective has sparse gradient almost everywhere, which produces valleys. These properties are considered unfavorable extreme cases for gradient methods, but are particularly common in tropical optimization problems. While the tropical linear regression problem is not convex, it is quasi-convex along tropical line segments; that is, its {sub-level} sets are tropically convex.

\subsubsection{Quasi-Convexity of Tropical Functions}

In formalizing tropical location problems, the work of \cite{comuaneci2024location} first considers functions which increase along tropical geodesics from a kernel $\vec v$; that is, for any point $\vec t$ on a geodesic from $\vec v$ to $\vec w$, we have that $f(\vec t) \leq f(\vec w)$. This property is equivalent to star-convexity of {sub-level} sets. To formally state this, we define (oriented) geodesic segments and $\triangle_{\min}$-star convexity of sets.

\begin{definition}[Oriented Geodesic Segment]
	The \emph{(oriented) geodesic segment} between $\vec a,\vec b \in \TPT{ N}$, is given by $[\vec a,\vec b]_{\triangle_{\min}} \coloneqq \{ \vec x \in \TPT{ N}: d_{\triangle_{\min}}(\vec a,\vec x)+ d_{\triangle_{\min}}(\vec x,\vec b) = d_{\triangle_{\min}}(\vec a,\vec b)\}$.
\end{definition}

The oriented geodesic segment is the union of all geodesics from $\vec a$ to $\vec b$ with respect to $d_{\triangle_{\min}}$.

\begin{definition}[$\triangle_{\min}$-Star Convex Set]
	A set $ S \subseteq \TPT{ N}$ is a \emph{$\triangle_{\min}$-star-convex} set with kernel $\vec v$ if, for every $\vec w \in  S$, we have $[\vec v,\vec w]_{\triangle_{\min}} \subseteq  S$.
\end{definition}

Using the above definitions, we can characterize a loss function which is increasing along tropical geodesics as a function with $\triangle_{\min}$-star-convex {sub-level} sets. The following definition gives a closed form expression for such functions, while \Cref{thm:sublevel_characterisation} proves their equivalence.

\begin{definition}[$\triangle_{\min}$-Star-Quasi-Convex Functions \citep{comuaneci2024location}]\label{def:star-quasi-convex}
	A function $f: \TPT{ N} \rightarrow \R$ is \emph{$\triangle_{\min}$-star-quasi-convex} with kernel $\vec v$ if $f(\vec x) = \hat{\gamma}(\vec x-\vec v)$ for some $\gamma: \R^{ N}_{\geq 0} \rightarrow \R$ which is increasing in every coordinate, where we define $\hat{\gamma}(\vec x) \coloneqq \gamma(\vec x - (\min_i x_i)\mathbf{1})$.
\end{definition}

\begin{theorem}[Theorem 17 of \cite{comuaneci2024location}]
	\label{thm:sublevel_characterisation}
	A continuous function $f: \TPT{ N} \rightarrow \R$ is $\triangle_{\min}$-star-quasi-convex with kernel $\vec v$ if and only if all of its non-empty {sub-level} sets are $\triangle_{\min}$-star convex with kernel $\vec v$.
\end{theorem}

\begin{remark}
	We can define $\triangle_{\max}$-star-quasi-convex functions similarly, but we also note that a function $f(\vec t)$ is a $\triangle_{\max}$-star-quasi-convex function if and only if $f(-\vec t)$ is a $\triangle_{\min}$-star-quasi-convex function.
\end{remark}

These $\triangle_{\min}$-star-convex functions {serve as a loss function with respect to a single data point}. For example, the tropical metric $\dtr$ is both $\triangle_{\min}$ and $\triangle_{\max}$-star-quasi-convex. The distance to a hyperplane $\dtr(\vec x, \mathcal{H}_{\vec t})$ is $\triangle_{\max}$-star-quasi-convex in $\vec x$, but $\triangle_{\min}$-star-quasi-convex in $\vec t$.

As a final remark, we note that $\triangle_{\min}$-star-quasi-convex functions are a special case of tropically quasi-convex functions by the following lemma. We limit our theoretical considerations to $\triangle_{\min}$-star-quasi-convex functions to utilize their closed form expression given in \Cref{def:star-quasi-convex}, however numerical experiments demonstrate strong performance of our methodology for more general tropical quasi-convex functions.

\begin{proposition}[Proposition 8 of  \cite{comuaneci2024location}]\label{prop:star_to_trop_convexity}
	Any $\triangle_{\min}$-star-convex set is $\min$-tropically convex. Hence, any $\triangle_{\min}$-star-quasi-convex function is \emph{tropically quasi-convex} in that it has tropically convex {sub-level} sets.
\end{proposition}

\subsubsection{Tropical Location Problems}

The $\triangle_{\min}$-star-convex functions defined above are generally a measure of closeness to the kernel $\vec v$, but in minimizing a statistical loss function we {consider proximity} to a sample of points. A \emph{location problem} \citep{laporte2019introduction} measures the closeness to a sample rather than a single kernel.

\begin{definition}[Min(/max)-Tropical Location Problem \citep{comuaneci2024location}]
	Consider a dataset $X = \{ \vec x_1, \dots, \vec x_K\} \subset \TPT{ N}$. For $k \leq K$, let $h_k$ be a $\triangle_{\min}$-star-quasi-convex function with kernel $\vec x_k$, and let $g: \R^K \rightarrow \R$ be increasing in every coordinate. A \emph{min-tropical location problem} is the minimization of an objective function $f: \TPT{ N} \rightarrow \R$ given by $f = g(h_1, \dots, h_K)$. A \emph{max-tropical location problem} is defined similarly, where the $h_k$ are $\triangle_{\max}$-star-quasi-convex.
\end{definition}

The main result of \cite{comuaneci2024location} is the following: tropical location problems contain a minimum in the tropical convex hull of its dataset $X$.

\begin{theorem}[Theorem 20 of \cite{comuaneci2024location}]\label{thm:min_in_tconv}
	Let $f$ be the objective function of a min-tropical location problem. Then there is a minimum of $f$ belonging to $\tconv_{\max}(X)$.
	Conversely, max-tropical location problems have minima belonging to $\tconv_{\min}(X)$.    
\end{theorem}

We can now revisit our motivating example (\Cref{def:lin_reg}), showing that tropical linear regression is not just a min-tropical location problem, but also demonstrates tropical quasi-convexity in that it has tropically convex {sub-level} sets.

\begin{proposition} \label{prop:LR_convexity}
	The tropical linear regression problem is a $\min$-tropical location problem with tropically convex {sub-level} sets.
\end{proposition}

{\it Proof}
\begin{align*}
	f(\vec t) &= {\max_{k \leq K} \left( \textstyle \max_i (\vec x_k-\vec t)_i - \text{2ndmax}_i (\vec x_k - \vec t)_i  \right)} \\
	&= \max_{k \leq K}\text{2ndmin}_j (\max_i (x_{ki} - t_i ) - x_{kj} + t_j ) \\
	&= \max_{k \leq K}\text{2ndmin}_j (t_j - x_{kj} - \min_i (t_i - x_{ki}{))}
\end{align*}
This is a min-tropical location problem with:
\begin{align*}
	g(\vec h) &= \max_{k \leq K} h_k, \\
	h_k(\vec t) &= \hat{\gamma}(\vec t - \vec x_k), \\
	\gamma(\vec t) &= \text{2ndmin}_{i}t_i.
\end{align*}
This $g$ and $\gamma$ are increasing in each coordinate and each $h_k$ has kernel $\vec x_k$, satisfying the conditions for a min-tropical location problem.

As the $h_k$ are $\triangle_{\min}$-star-convex, their {sub-level} sets are min-tropically convex by \Cref{prop:star_to_trop_convexity}. This is preserved when taking a maximum over samples, so $f$ also has $\min$-tropically convex {sub-level} sets.
\qed \\

As mentioned in \Cref{prop:star_to_trop_convexity}, $\triangle_{\min}$-star-quasi-convexity implies tropical quasi-convexity. However, the tropical linear regression problem is a case in which the reverse implication does not hold; in general, the linear regression loss function is not $\triangle_{\min}$-star-quasi-convex.

\subsection{Tropical Data Science} \label{subsec:trop_data_science}

While young, the field of tropical data science is motivated by the crucial question of statistical learning from the increasing volume of phylogenetic tree data \citep{kahn2011future}. In contrast to the BHV \citep{billera2001geometry} or Robinson--Foulds metrics \citep{ROBINSON1981131}, the tropical interpretation of tree space provides both interpretable geometry and computational efficiency \citep{penny1993distributions,monod2018StatPerspective}. Many of the statistical methods defined for the tropical setting are framed as {optimization} problems \citep{barnhill2023tml}, and in this subsection we review some such tropical statistical optimization problems which we use to test our tropical descent methodology. We note the varying degrees of convexity demonstrated by each of them. While we restrict our considerations through the rest of this paper to statistics which arise from tropical location problems, in \Cref{appsec:further_experiments} we discuss other tropical statistical optimization problems.

\subsubsection{Centrality Statistics}

Our first problems of interest, Fermat--Weber points and Fr\'echet means, are centrality statistics for data on a general metric space; Fermat--Weber points \citep{wesolowsky1993weber} act as a generalized median, while Fr\'echet means \citep{frechet1948elements} are a generalization of the mean.

\begin{definition}[Tropical Fermat--Weber Points \citep{lin2018FermatWeber}]
	Let $X = \{\vec x_1, \dots, \vec x_K \} \subset \TPT{N}$ be our dataset. A solution to the following optimization problem is a \emph{tropical Fermat--Weber point}:
	\begin{align*}
		\min_{\vec t \in \TPT{N}} f(\vec t) = \frac{1}{K}\sum_{k \leq K} d_{\tr}(\vec x_k, \vec t).
	\end{align*}
\end{definition}

\begin{definition}[Tropical Fr\'echet means \citep{monod2018StatPerspective}]
	Let $X = \{\vec x_1, \dots, \vec x_K \} \subset \TPT{N}$ be our dataset. A solution to the following optimization problem is a \emph{tropical Fr\'echet mean}:
	\begin{align*}
		\min_{\vec t \in \TPT{N}} f(\vec t) = \left[\frac{1}{K}\sum_{k \leq K} d_{\tr}(\vec x_k, \vec t)^2 \right]^{1/2}.
	\end{align*}
\end{definition}

The Fermat--Weber and Fr\'echet mean problems are both max-tropical and min-tropical location problems, as the tropical metric $\dtr$ is both $\triangle_{\min}$ and $\triangle_{\max}$-star-quasi-convex. In fact, they are also both classically convex and hence act as a benchmark against which we can test our gradient methods.

\subsubsection{Wasserstein Projections}

The next statistical optimization problem was motivated by \cite{cai2022distances}, and looks to identify a tropical projection which minimizes the Wasserstein distance between samples in different tropical projective tori.

\begin{definition}[Tropical Wasserstein Projections \citep{talbut2023probability}]
	Let $X = \{ \vec x_1, \dots, \vec x_K \} \subset \TPT{N}$ be a dataset in one space, and $Y = \{ \vec y_1, \dots, \vec y_K \} \subset \TPT{M}$ be a dataset in a different space. Let $J_1, \dots, J_M$ be a non-empty partition of $[N]$. Then the \emph{tropical $p$-Wasserstein projection problem} is the minimization of the following objective function:
	\begin{align*}
		\min_{\vec t \in \TPT{N}}  f_p(\vec t) = \left(\frac{1}{K}\sum_{k \leq K} \|  (\max_{i \in J_j} (\vec x_k - \vec t)_i )_{ {j \leq M}} - \vec y_k\|_{\mathrm{tr}}^p\right)^{1/p}
	\end{align*}
	The \emph{tropical $\infty$-Wasserstein projection problem} is the minimization of:
	\begin{align*}
		\min_{\vec t \in \TPT{N}}  f_{\infty}(\vec t) = \max_{k \leq K} \|  ( \max_{i \in J_j} (\vec x_k - \vec t)_i  )_{ {j \leq M}} - \vec y_k\|_{\mathrm{tr}}
	\end{align*}
\end{definition}
In contrast to the work by \cite{lee2021tropical} which looks to solve the optimal transport problem for intrinsic measures on $\TPT{N}$, this $p$-Wasserstein projection problem is computing an optimal projection for empirical samples.

\begin{proposition}
	The tropical $p,\infty$-Wasserstein projections are min-tropical location problems, and the $\infty$-Wasserstein projection has $\min$-tropically convex {sub-level} sets.
\end{proposition}

{\it Proof}
Let $\vec z_k = (x_{ki} - y_{kj})_{i \leq N}$ where $j$ is such that $J_j$ is the unique partition set containing $i$. We then have that 
\begin{align*}
	\|  (\max_{i \in J_j} (\vec x_k - \vec t)_i )_{ {j \leq M}} - \vec y_k\|_{\mathrm{tr}} 
	&=\max_{b} \max_{a \in J_b}  ( z_{ka} - t_a  ) - \min_{j} \max_{i \in J_j}  (z_{ki} - t_i ) \\
	&= \max_{a}  (z_{ka} - t_a  ) + \max_j\min_{i \in J_j} ( t_i - z_{ki}  ) \\
	&= \max_j  (\min_{i \in J_j}  (t_i - z_{ki} - \min_{a}  (t_a - z_{ka} {)))}
\end{align*}
The $p, \infty$-Wasserstein projections can then be written as
\begin{align*}
	f_{p} (\vec t) &=\left(\frac{1}{K}\sum_{k \leq K} \left( \max_j  (\min_{i \in J_j} (t_i - z_{ki} - \min_{a}  (t_a - z_{ka} {)))} \right)^p\right)^{1/p} \\
	f_{\infty} (\vec t) &=\max_{k \leq K} \max_j  (\min_{i \in J_j}  (t_i - z_{ki} - \min_{a} (t_a - z_{ka}{)))}
\end{align*}
These are max-tropical location problems with:
\begin{align*}
	g_{p}(\vec h) &= \left(\frac{1}{K}\sum_{k \leq K} h_k^p\right)^{1/p}, \\
	g_{\infty}(\vec h) &= \max_{k \leq K} h_k, \\
	h_k(\vec t) &= \hat{\gamma}(\vec t - \vec z_k), \\
	\gamma(\vec t) &= \max_j\min_{i \in J_j}t_i.
\end{align*}
The $g_p, g_{\infty}$ and $\gamma$ are increasing in each coordinate and each $h_k$ has kernel $\vec z_k$, satisfying the conditions for a min-tropical location problem.

As the $h_k$ are $\triangle_{\min}$-star-convex, their {sub-level} sets are min-tropically convex by \Cref{prop:star_to_trop_convexity}. This is preserved when taking a maximum over samples, so $f_{\infty}$ also has $\min$-tropically convex {sub-level} sets.
\qed \\

\begin{figure}
	\centering
	\includegraphics[width = 0.5\textwidth]{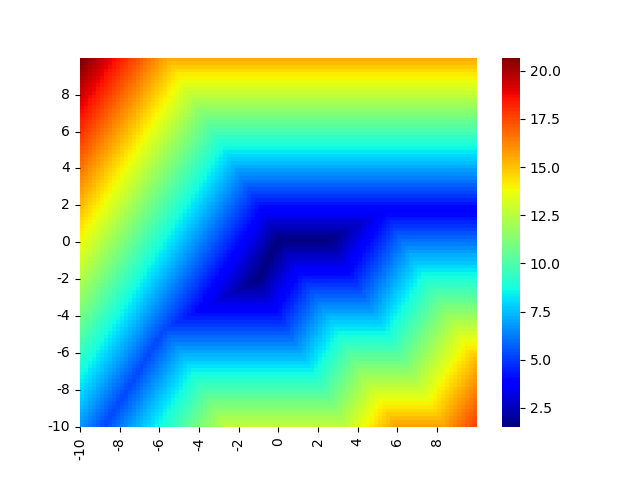}
	\caption{A heat map of a 3-sample $\infty$-Wasserstein projection objective on $\TPT{3} \cong \mathcal{H} = \{ \vec x: \sum x_i = 0 \}$.}
	\label{fig:WDDDinf_Visualisation}
\end{figure}

\Cref{fig:WDDDinf_Visualisation} shows a heatmap for the $\infty$-Wasserstein projection objectives. As in the tropical linear regression problem, it is piecewise linear, non-convex and has sparse gradient.

\section{Steepest Descent with the Tropical Norm}\label{sec:contribution}

This section comprises the theoretical contributions of this paper; we outline our proposed method of tropical descent, and the theoretical convergence guarantees for tropical descent when applied to tropical location problems.

We treat optimization over $\TPT{N}$ as optimization over $\R^{N}$ with the knowledge that $f$ is constant along $\mathbf{1}$-rays. We assume we can compute some gradient-like function $\nabla f: \R^N \rightarrow \R^N$ which {satisfies} \Cref{ass:derivatives}.

\begin{assumption}\label{ass:derivatives}
	Let $f = g(h_1, \dots, h_K)$ be a min-tropical location problem. Then we have that
	\[
	\nabla f(\vec t)_j < 0 \Rightarrow \exists \, k \text{ s.t. } j \in \argmin_i(t_i - x_{ki}).
	\]
\end{assumption}

This is a very general assumption, which holds for some of the most natural gradient-like functions. For example, if $\nabla f(\vec t)$ is given by some Clarke gradient, that is, some element of the \emph{generalized Clarke derivative} of $f$, then \Cref{lemma:quasi-convex_derivatives} shows that \Cref{ass:derivatives} holds. 

\begin{definition}[Generalized Clarke Derivative (2.13, I.I of \cite{kiwiel2006methods})]
	We define the set $M_f(\vec t)$ by
	\[
	M_f(\vec t) = \{ \vec z \in \R^N : \nabla f(\vec y_l) \rightarrow \vec z, \text{ where }\vec y_l \rightarrow \vec t, f\text{ differentiable at each }\vec y_l \}.
	\]
	Then the \emph{generalized Clarke derivative} is given by
	\[
	\partial f(\vec t) = \conv M_f(\vec t).
	\]
\end{definition}

\begin{lemma}\label{lemma:quasi-convex_derivatives}
	Suppose $f = g(h_1, \dots, h_K)$ is min-tropical location problem, and for all $\vec t$, $\nabla f(\vec t)$ is some element of $\partial f(\vec t)$. Then \Cref{ass:derivatives} holds.
\end{lemma}

{\it Proof}
Suppose that $\nabla f(\vec t)_j < 0$. Then there is a sequence of differentiable $\vec y_\ell$ satisfying $\nabla f(\vec y_\ell)_j < 0$, $\vec y_\ell \rightarrow \vec t$. Hence by the chain rule, as $g$ is increasing in every $h_k$, there is some $h_k$ which is strictly decreasing in the $j^{th}$ coordinate at $\vec y_\ell$. WLOG, the same $h_k$ is strictly decreasing in the $j^{th}$ coordinate for all $\vec y_\ell$. As $\gamma$ is increasing in every coordinate, $h_k = \hat{\gamma}(\vec y_\ell - \vec x_k)$ is strictly decreasing in the $j^{th}$ coordinate if and only if $j \in \argmin_i (y_{\ell i} - x_{ki})$. This is preserved through the limit as $\ell \rightarrow \infty$, so we conclude that $j \in \argmin_i  (t_i - x_i  )$.
\qed \\

In practice, rather than computing elements of the {generalized} Clarke derivative, we use PyTorch automatic differentiation \citep{paszke2019pytorch}. The derivatives computed by PyTorch are dependent on the computational representation of a function, but are designed to respect the chain rule. We must therefore ensure that our computational representation of $g$ has non-negative derivatives in each $h_k$, while each $h_k$ is computed such that its $t_j$ derivative is negative only if $j \in \argmin_i (t_i - x_{ki})$. For such a computational representation of $f$, \Cref{ass:derivatives} will hold by the chain rule. For all the tropical location problems we consider in this paper, we find their most natural computational representation gives rise to derivatives $\nabla f$ satisfying \Cref{ass:derivatives}.

\subsection{Proposed Method}\label{subsec:theory_idea}

In defining tropical descent, we refer back to the motivation behind classical gradient descent---steepest descent directions.

\begin{definition}[Steepest Descent Direction \citep{boyd2004convex}]
	A normalized \emph{steepest descent direction} at $\vec t$ with respect to the norm $\| \cdot \|$ is:
	\[
	\vec d = \argmin \{ \nabla f(\vec t)^{\top} \vec v : \|\vec v\| = 1 \}.
	\]
\end{definition}

With respect to the Euclidean norm, the normalized derivative is a steepest descent direction. {The steepest descent directions under the tropical norm are characterized by} the following lemma.

\begin{lemma}\label{lemma:Trop_Descent_Directions}
	Suppose $\nabla f \neq \vec 0$. Then the direction $\vec v$ is a normalized tropical steepest descent direction {(up to $\sim$ equivalence)} iff
	\begin{align*}
		\nabla f_i < 0 & \Rightarrow v_i = 1, \\
		\nabla f_i > 0 & \Rightarrow v_i = 0.
	\end{align*}
\end{lemma}

{\it Proof}
WLOG, we assume $0 \leq v_i \leq 1$. Then
\[
\nabla f^{\top} \vec v \geq \sum_{i:\,\nabla f_i < 0} \nabla f_i
\]
This is achieved if and only if $v_i = 1$ for all $i$ such that $\nabla f_i <0$, and $v_i = 0$ for all $i$ such that $\nabla f_i >0$. We note that any such $v_i$ satisfies $\|\vec v\|_{\tr} = 1$ as $\textstyle \sum_i \nabla f_i = 0$ by the equivalence relation of the tropical projective torus, and $\nabla f_i \neq \vec 0$ so $\nabla f$ has positive and negative coordinates.
\qed \\

\Cref{lemma:Trop_Descent_Directions} uniquely defines a descent direction $\vec d$ if and only if, for every coordinate $i$, $\nabla f_i \neq 0$. To uniquely define a descent direction for any $\vec t$, we take the direction $\vec d$ given by
\begin{align*}
	d_i = \begin{cases}
		1 & \text{ if }\nabla f_i < 0, \\
		0 & \text{ otherwise.}
	\end{cases}
\end{align*}

We refer to this as a \emph{min-tropical descent} direction. In contrast, 
\begin{align*}
	d_i = \begin{cases}
		-1 & \text{ if }\nabla f_i > 0, \\
		0 & \text{ otherwise.}
	\end{cases}
\end{align*}
defines a \emph{max-tropical descent} direction.

Throughout the rest of this paper, we use the min-tropical descent direction as convention, and will specify when the max-tropical descent direction is used instead. Tropical descent therefore refers to the method outlined by \Cref{alg:TropDescent}, for some specified step size sequence $(a_{ m})_{ m\geq 1}$, while we refer to classical gradient descent as classical descent.

\begin{algorithm}[ht!]
	\caption{Tropical Descent}\label{alg:TropDescent}
	\begin{algorithmic}
		\For{$ m < $ max\_steps}
		
		\State $\vec g_{ mi} \gets \nabla f(\vec t_{ m-1})_i$;
		\State $\vec d_{ mi} \gets \mathbf{1}_{\vec g_{ mi} < 0}$;
		\State $\vec t_{ m} \gets \vec t_{ m-1} + a_{ m} \vec d_{ m}$;
		\EndFor
	\end{algorithmic}
\end{algorithm}

We note that over a large number of steps, \Cref{alg:TropDescent} can produce large $t_{ mi}$ values as each coordinate is non-decreasing in $ m$. If this results in overflow or precision errors, a {normalization} operation can be included in each step, such as subtracting $\frac{1}{N} \sum_{{i \leq N} } t_{ mi}$ from each $t_{ mi}$. The function $f$ is invariant under such a translation due to the equivalence relation on the tropical projective torus.

We note that \cite{barnhill2024polytropes} and \cite{barnhill2023tml} implement classical descent for the computation of Fermat--Weber points; this is the most natural choice as the Fermat--Weber problem is classically convex. However, the linear regression and Wasserstein projection problems lack classical convexity, and hence do not consistently produce accurate solutions via classical descent. It is these tropically convex problems for which tropical descent is necessary.

\subsection{Convergence Guarantees}\label{subsec:convergence}

In this section we prove the theoretical guarantees of tropical descent for min-tropical location problems. We show that tropical descent must converge to the tropical convex hull of the data, and hence will necessarily find a minimum of a $\triangle_{\min}$-star-convex problem. Throughout this section, we assume the step sizes $a_{ m}$ satisfy $a_{ m} \rightarrow 0, \sum a_{ m} = \infty$. 

We first highlight the key difference in stability between classical descent and tropical descent; if some component of the derivative is locally zero, tropical descent cannot converge in that neighborhood.

\begin{proposition} \label{prop:instability}
	Let $\vec t_{ m}$ be some sequence defined by tropical descent. Assume that for all $ m$, $\nabla f(\vec t_{ m}) \neq \vec 0$. Suppose in some bounded open $U$, $\nabla f_i$ is identically 0. Then no point in $U$ is a convergence point for $\vec t_{ m}$; in particular, if $\vec t_m \in U$ there is some $n > m$ such that $\vec t_n \notin U$.
\end{proposition}

{\it Proof}
Assume otherwise. Then for all $n \geq m$, we have $\nabla f(\vec t_n)_i = 0$, so $t_{ni}$ is constant. At each step, we increase at least one coordinate by $a_n$, so as $n \rightarrow \infty$:
\[
\sum_j t_{nj} \geq \sum_{j} t_{mj} + \sum_{n \geq m' \geq m} a_{m'} \rightarrow \infty.
\]
Therefore $\|\vec t_n\|_{\tr}$ will become arbitrarily large, as $\max_j t_{nj} \rightarrow \infty$ while $\min_j t_{nj} \leq  t_{mi}$. Then we must have left $U$; contradiction.
\qed \\

As a result of the proposition above, local minima such as the valleys in \Cref{fig:LR_Visualisation,fig:WDDDinf_Visualisation} can be stable with respect to classical descent but unstable with respect to tropical descent.

We next prove the main convergence result of this work; tropical descent must converge to the tropical convex hull of the data points.

\begin{theorem}\label{thm:tconv_convergence}
	Suppose $f$ is a min-tropical location problem with respect to the dataset $X = \{\vec x_1, \dots, \vec x_K \}$. Let $(\vec t_{ m})_{ m \geq 1}$ be a sequence of points defined by tropical descent, such that for all $ m$: $\nabla f(\vec t_{ m}) \neq 0$. Let $V = \tconv_{\max}(\vec x_1, \dots, \vec x_K)$, and define
	\begin{align*}
		\Delta_j(\vec t) = \sum_{ {k \leq K}} \sum_{{i \leq N}} [t_j - x_{kj} - t_i + x_{ki}]^+.
	\end{align*}
	For any $\epsilon > 0$, let $M_1$ be such that for all $m \geq M_1$, $a_m < \epsilon/N$. Let $M_2$ be such that for all $m \geq M_2$ $s_{m} \geq \max_j \Delta_j(\vec t_{M_1}) + s_{M_1} - K\epsilon$, where $s_m$ is the partial sum of the step size sequence $a_m$. Then for all $m \geq M_2$ we have:
	\[
	d_{\tr}(\vec t_{ m}, V) \leq \epsilon.
	\]
\end{theorem}

{\it Proof}
For each $i \leq N$, we define the functions
\[
\delta_j(\vec t) = \min_{ {k \leq K}} \left[\sum_{ {i \leq N}} [t_j - x_{kj} - t_i + x_{ki}]^+ \right].
\]
\begin{claim}\label{claim:delta_vs_V}
	For all $\vec t \in \TPT{N}$, we have $d_{\tr}(\vec t,V) \leq \max_j \delta_j(\vec t)$.
\end{claim}
{\it Proof}
As discussed by \cite{maclagan2015introduction}, there is a well-defined metric projection map $\pi_V$ onto $V$ given by:
\begin{gather*}
	\pi_V(\vec t) = \lambda_1 \odot \vec x_1 \boxplus \dots \boxplus \lambda_K \odot \vec x_K, \\
	\text{ where }\lambda_k = \max \{ \lambda \in \R : \, \lambda \odot \vec x_{k} \boxplus \vec t = \vec t \} = \min_i(t_i - x_{ki}).
\end{gather*}
The distance from $\vec t$ to $V$ is then given by $d(\vec t, \pi_V(\vec t))$, and we have:
\begin{align*}
	d(\vec t, V) &= d(\vec t, \pi_V(\vec t)) \\
	&= \| (\max_k (\min_i (t_i - x_{ki}) - t_j + x_{kj})_{j \leq N} \|_{ \tr}.
\end{align*}
We note that for all $j,k$, the term $\min_i (t_i - x_{ki}) - t_j + x_{kj}$ is negative. Hence:
\begin{align*}
	d(\vec t, V) & \leq 0 - \min_j (\max_k (\min_i (t_i - x_{ki}) - t_j + x_{kj})) \\
	&  = \max_j (\min_k (t_j - x_{kj} - \min_i (t_i - x_{ki}))) \\
	& \leq \max_j \delta_j(\vec t).
\end{align*}
The claim is proved.
\qed \\
We will now prove that for all $m \geq M_2$, we have that $\delta_j(\vec t_m) \leq \epsilon$. 
\begin{claim} \label{claim:delta_jumps}
	If $\delta_i(\vec t_m) = 0$ and $m \geq M_1$, then $\delta_i(\vec t_{m+1}) < \epsilon$.
\end{claim}

{\it Proof}
For all ${k \leq K, j \leq N}$, we have that $t_{mj} - x_{kj} - t_{mi} + x_{ki}$ is $1$-Lipschitz in $\vec t$ (with respect to $\dtr$). Therefore $\delta_j$ is $N$-Lipschitz in $\vec t$. By the tropical descent rule:
\begin{align*}
	\delta_j(\vec t_{m+1}) &\leq \delta_j(\vec t_m) + N\dtr(\vec t_{m}, \vec t_{m+1})\\
	&= 0 + Na_m \\
	&< \epsilon.
\end{align*}
The claim is proved.
\qed \\

\begin{claim}\label{claim:delta_decreases}
	If $\delta_j(\vec t_m) > 0$ then $\delta_j(\vec t_{m+1}) \leq \delta_j(\vec t_{m})$.
\end{claim}

{\it Proof}
Suppose $\delta_j(\vec t_m) > 0$. Then $\forall k$, $t_{mj} - x_{kj} > \min_i (t_{mi} - x_{ki})$, so by \Cref{lemma:quasi-convex_derivatives} we have $\nabla f_j \geq 0$. The tropical steepest descent direction therefore has $d_{mj} = 0$. We conclude that
\begin{align*}
	\forall \, k,i : \quad &t_{(m+1)j} - x_{kj} - t_{(m+1)i} + x_{ki} \leq t_{mj} - x_{kj} - t_{mi} + x_{ki}, \\
	&\Rightarrow \quad \delta_j(\vec t_{m+1}) \leq \delta_j(\vec t_m).
\end{align*}
The claim is proved.
\qed \\

\begin{claim}\label{claim:delta_to_zero}
	For some $M_1 \leq m \leq M_2$,  $\delta_j(\vec t_m) < \epsilon$.
\end{claim}

{\it Proof}
We note that $K\delta_j(\vec t) \leq \Delta_j(\vec t)$; hence if $\Delta_j(\vec t_{M_1}) < K\epsilon$ then we are done.

We assume $\delta_j(\vec t_m) > \epsilon$ for all $M_1 \leq m \leq M_2$, and prove the result by contradiction. Then as in the proof of the previous claim, at each step we have $\nabla f_j \geq 0$, $d_{mj} = 0$ and $t_{mj}$ is fixed for $M_1 \leq m \leq M_2$.

We have assumed that for each $m$, $\nabla f \neq \vec 0$ and so there is some $i$ such that $\nabla f_i < 0$, and by \Cref{lemma:quasi-convex_derivatives}, there is some $k$ such that $i \in \argmin_{\ell}  ( t_{m\ell} - x_{k\ell} )$. Hence:
\begin{align*}
	\forall \, \ell: \quad t_{mj} - x_{kj} - t_{mi} + x_{ki} &\geq t_{mj} - x_{kj} - t_{m\ell} + x_{k\ell}, \\
	\Rightarrow \quad t_{mj} - x_{kj} - t_{mi} + x_{ki} &\geq \delta_j(\vec t_m)/ N, \\
	&> \epsilon/N, \\
	&\geq a_m.
\end{align*}
Hence this term in the summation of $\Delta_j(\vec t_m)$ decreases by $a_m$, while other terms are non-increasing as $t_{mj}$ is fixed. 

Noting that $K\delta_j(\vec t) \leq \Delta_j(\vec t)$, we conclude
\begin{align*}
	K\delta_j(\vec t_{M_2}) &\leq \Delta_j (\vec t_{M_2}), \\
	&\leq \Delta_j (\vec t_{M_1}) - s_{M_2} + s_{M_1}, \\
	&\leq K\epsilon.
\end{align*}
Contradiction. So the claim is proved.
\qed \\

Finally, from the claims above, we conclude that $ \forall m \geq M_2, {i \leq N}$, we have $\delta_i(\vec t_m) \leq \epsilon$, and so by our first claim, we have $d_{\tr}(\vec t_m, V) \leq \epsilon$.
\qed \\

As proven by \cite{comuaneci2024location}, tropical location problems have minima in the tropical convex hull of the data. The theorem above says that according to tropical descent, only such minima can be stable. We cannot do much better in terms of convergence results for general location problems; they can have disconnected {sub-level} sets so we would not expect gradient methods to find a global minimum. We would require some degree of convexity to guarantee global solvability, which is exactly what the following corollary gives us; tropical descent will necessarily converge to a minimum of any $\triangle_{\min}$-star-quasi-convex function.

\begin{corollary}\label{cor:1sample_convergence}
	Let $f(\vec t)$ be a $\triangle_{\min}$-star-quasi-convex function, and let $(\vec t_{ m})_{ m \geq 1}$ be a sequence of points defined by tropical descent such that for all $ m$: $\nabla f(\vec t_{ m}) \neq 0$. Then there is some global minimum $\vec t^*$ such that for all $m  \geq M_2$, $d_{\tr}(\vec t_m, \vec t^*) \leq \epsilon$ where $M_2$ is as in \Cref{thm:tconv_convergence}.
\end{corollary}

{This corollary follows from considering the 1-sample case of tropical location problems; \cref{cor:1sample_convergence} establishes convergence to the kernel of a $\triangle_{\min}$-star-quasi-convex function, which must necessarily be a global minimum. This convergence result holds even if $f$ does not have a unique minimum --- it may be that $f$ has multiple global minima but a unique kernel, or $f$ may have a full dimensional set of global minima upon which $\nabla f(\vec t) = \vec 0$ and tropical descent terminates.}\\

As a further corollary, we note that we can achieve the same result for max-tropical location problems by using max-tropical descent directions.

\begin{corollary}\label{cor:max-trop_problem}
	Suppose $f$ is a max-tropical location problem with respect to the dataset $X = \{\vec x_1, \dots, \vec x_K \}$. Let $(\vec t_{ m})_{ m \geq 1}$ be a sequence of points defined via max-tropical descent directions, such that for all $ m$: $\nabla f(\vec t_{ m}) \neq 0$. Let $V = \tconv_{\min}(\vec x_1, \dots, \vec x_K)$, and define
	\[
	\Delta'_i(\vec t) = \sum_{{k \leq K}} \sum_{{j \leq N}}  [t_j - x_{kj} - t_i + x_{ki}]^+.
	\]
	For any initial $\vec t_0$ and $\epsilon > 0$, let $M_1$ be such that for all $m \geq M_1$, $a_m \leq \epsilon/N$. Let $M_2$ be such that for all $m \geq M_2$ $s_{m} \geq \max_i \Delta'_j(\vec t_{M_1}) + s_{M_1} - K\epsilon$. Then for all $m \geq M_2$ we have:
	\[
	d_{\tr}(\vec t_{ m}, V) \leq \epsilon.
	\]
\end{corollary}

\begin{remark}
	The results above indicate convergence to $V$ or $\vec t^*$ in the tropical topology of $\TPT{N}$, but $\vec t_m$ cannot converge in the Euclidean topology of $\R^N$ due to the divergent partial sums of $a_m$. 
\end{remark}

We conclude by establishing {a foundational bound for} the convergence rate of tropical descent; we take step sizes given by $a_m = \alpha m^{-1/2}\|\nabla f(\vec t_m) \|_{\tr}$, and assume derivatives have bounded tropical norm as is the case for piecewise linear 1-Lipschitz objective functions {(with finitely many linear pieces)}. In this case, we prove convergence at a rate of  {$O(1/\sqrt{m})$}.

\begin{proposition}\label{prop:convergence_rate}
	Suppose $f$ is a min-tropical location problem as in \Cref{thm:tconv_convergence}. Assume that there is some $L \in \N$ such that for all $\vec t_m$, $1/L \leq \|\nabla f(\vec t_m)\|_{\tr} \leq 2$, and assume further there is some $D$ such that {$\forall m \in \N: \, \sum_{k \leq K} d_{\triangle_{\max}}(\vec x_k, \vec t_{m}) \leq KD$}. We define $a_m = \alpha m^{-1/2}\|\nabla f(\vec t_m) \|_{\tr}$.
	Then for all $m$ satisfying
	\begin{align} \label{eq:minimum_m}
		{ m \geq \max \left\{ 
			\left(\sqrt 2 + \frac{KDL}{2\alpha} - \sqrt{LKN-1} \right)^2-1, 
			\left(\frac{2 \alpha N}{D}\right)^2 
			\right\}, }
	\end{align}
	we have that:
	\[
	d_{\tr}(\vec t_m, V) \leq {\frac{2 \sqrt2 \alpha N}{\sqrt{m+1} - KDL/2\alpha + \sqrt{LKN-1} - \sqrt2} } = \epsilon_m.
	\]
	In particular, if $f$ is a $\triangle_{\min}$-star-quasi-convex function, then there is some global minimum $\vec t^*$ such that 
	\[
	d_{\tr}(\vec t_m, \vec t^*) \leq {\frac{2 \sqrt2 \alpha N}{\sqrt{m+1} - DL/2\alpha + \sqrt{LN-1} - \sqrt2}  = \epsilon'_m.}
	\]
\end{proposition}

{In order to use \cref{thm:tconv_convergence}, we will need to find some $M_1$ such that $a_{M_1} < \epsilon/N$ (our step size is small enough) and $s_{m} \geq \max_j \Delta_j(\mathbf t_{M_1}) + s_{M_1} - K\epsilon$ (we have taken enough steps to get close to $V$). Our condition that $m$ is large enough such that $\sqrt{m+1} - KDL/2\alpha + \sqrt{LKN-1} \geq \sqrt 2$ ensures that such an $M_1$ exists, while $\sqrt{m} \geq 2 \alpha N/D$ ensures that $M_1 \leq m$. We can then take the largest such $M_1$ to ensure that the corresponding $\epsilon$ is as tight as possible.}\\ 

{\it Proof}
We fix $m$ {satisfying \cref{eq:minimum_m}} and define: 
\begin{align*}
	{\xi(m) }&{\coloneqq \frac{\sqrt{m+1}}{2} - \frac{KDL}{4\alpha} + \sqrt{\left(\frac{\sqrt{m+1}}{2} - \frac{KDL}{4\alpha}\right)^2 + LKN-1},} \\
	M_1 &\coloneqq \left \lfloor \xi(m)^2 \right \rfloor.
\end{align*}
{ We first note that by the quadratic mean inequality and as $\sqrt{m+1} - KDL/2\alpha + \sqrt{LKN-1} \geq \sqrt 2$:
	\begin{align*}
		\xi(m) &= \frac{\sqrt{m+1}}{2} - \frac{KDL}{4\alpha} + \sqrt{\left(\frac{\sqrt{m+1}}{2} - \frac{KDL}{4\alpha}\right)^2 + LKN-1} \\
		&\geq \frac{\sqrt{m+1}}{2} - \frac{KDL}{4\alpha} + \frac{1}{\sqrt{2}} \left( \frac{\sqrt{m+1}}{2} - \frac{KDL}{4\alpha} + \sqrt{LKN-1} \right) \\
		&\geq \frac{1}{\sqrt 2} \left( \sqrt{m+1} - \frac{KDL}{2\alpha} +  \sqrt{LKN-1} \right) \geq 1.
	\end{align*}
}
Therefore we have that
\begin{gather*}
	{1} \leq \sqrt{M_1} \leq \xi(m),
\end{gather*}
and so $\sqrt{M_1}$ lies between the roots of the polynomial
\[
{ p(t) = L^{-1}t^2 + \left( KD/2\alpha - L^{-1}\sqrt{m+1}\right) t+ L^{-1} - KN}.
\]
Therefore:
\begin{gather*}
	{L^{-1}}(M_1+1) + \left( {KD /2\alpha} - {L^{-1}} \sqrt{m+1}\right) \sqrt{M_1}- KN \leq 0, \\
	{L^{-1}}\sqrt{M_1}\sqrt{M_1+1} + {KD\sqrt{M_1}/2\alpha} - {L^{-1}}\sqrt{M_1}\sqrt{m+1} \leq KN, \\
	{KD} + 2\alpha L^{-1}(\sqrt{M_1+1} - \sqrt{m+1}) \leq \frac{2\alpha KN}{\sqrt{M_1}}.
\end{gather*}
{
	We now show that $M_1 \leq m$. Suppose otherwise; then by the inequality above:
	\begin{align*}
		KD &\leq \frac{2 \alpha K N}{\sqrt{M_1}} + 2\alpha L^{-1}(\sqrt{m+1} - \sqrt{M_1+1}) \\
		&<  \frac{2 \alpha K N}{\sqrt{m}}
	\end{align*}
	This contradicts our assumption that $m \geq (2\alpha N/D)^2$, so we conclude that $M_1 \leq m$.}
\begin{claim}
	For all $n \geq {m}$ we have $a_n \leq \epsilon_m / N$.
\end{claim}
{\it Proof}
{We have shown that $\xi(m) \geq 1$, so:}
\begin{align*}
	M_1 &\geq \xi(m)^2 - 1, \\
	&\geq \left(\xi(m) -1\right)^2, \\
	\sqrt{M_1} &\geq \xi(m) - 1, \\
	&\geq {\frac{1}{\sqrt{2}} \left( \sqrt{m+1} - \frac{KDL}{2\alpha} + \sqrt{LKN-1} \right) - 1}.
\end{align*}
Hence, by the step size definition { and the fact that $n \geq m \geq M_1$}:
\begin{align*}
	a_n &\leq \frac{2 \alpha}{\sqrt{M_1}}, \\
	&\leq {\frac{2 \sqrt2 \alpha}{\sqrt{m+1} - KDL/2\alpha + \sqrt{LKN-1} - \sqrt2}} = \epsilon_m/N.            
\end{align*}
The claim is proved.
\qed
\begin{claim}
	We have $\max_j \Delta_j(\vec t_{M_1}) + s_{M_1} - s_m \leq K\epsilon_m $.
\end{claim}
{\it Proof}
We note that
{\begin{align*}
		\Delta_j(\vec t_{M_1}) &= \sum_{k \leq K} \sum_{i \leq N} [t_{M_1j} -x_{kj} - t_{M_1i} + x_{ki}]^+,\\
		&\leq \sum_{k \leq K} d_{\triangle_{\max}}(\vec x_k, \vec t_{M_1}),\\
		&\leq KD.
	\end{align*}
	And therefore we have $\max_j \Delta_j(\vec t_{M_1}) \leq KD$. }
The partial sums are bounded by:
\begin{align*}
	s_{M_1} - s_m &\leq \sum_{n=M_1+1}^m -a_n \\
	&\leq \sum_{n=M_1+1}^m \frac{-\alpha}{Ln^{1/2}} \\
	&\leq \sum_{n=M_1+1}^m - 2 \alpha L^{-1} (\sqrt{n+1} - \sqrt{n}) \\
	&\leq 2\alpha L^{-1} {(\sqrt{M_1+1} - \sqrt{m+1})}.
\end{align*}
Therefore:
\begin{align*}
	\max_i\Delta_i(\vec t_{M_1}) + s_{M_1} - s_m &\leq {KD} +2\alpha L^{-1}(\sqrt{M_1+1} - \sqrt{m+1}), \\
	&\leq \frac{2\alpha KN}{\sqrt{M_1}} \leq K\epsilon_m.
\end{align*}
\qed

Hence by \Cref{thm:tconv_convergence}, we have that 
\[
d_{\tr}(\vec t_m, V) \leq \epsilon_m.
\]
The bound for the tropically quasi-convex case follows from \Cref{cor:1sample_convergence}.
\qed \\

{Our error bound in Proposition 12 is a useful foundational result on the convergence behavior of tropical descent; for sufficiently nice problems, tropical descent will converge at a rate of $O(1/\sqrt{m})$. However, a more detailed study of the numerical analysis of tropical descent would be valuable. We outline possible questions below, as well as presenting a comparison of our error bounds with tropical descent for tropical linear regression.}

{Firstly, the constant term $LD/2\alpha$ in the denominator of $\epsilon_m$ will typically be larger than necessary, resulting in a weaker error bound. In particular, this constant term is quick to dominate for larger problems. Further work improving these constant factors would produce much tighter and more robust error bounds. Secondly, the error $\epsilon_m$ bounds the distance to a tropical convex hull $V$, which is rarely a useful quantity to bound. We would rather use $\epsilon'_m$, which bounds the distance to a global minimum, but our tropical location problems are not known to be $\triangle_{\min}$-star-quasi-convex. It would be particularly useful to generalize our error bound $\epsilon_m'$ to such problems which exhibit some local (but not global) $\triangle_{\min}$-star-quasi-convexity.}

{For a qualitative view on \cref{prop:convergence_rate}, we compare $\epsilon_m'$ and the convergence rate of tropical descent for tropical linear regression. As noted, tropical linear regression is tropically quasi-convex (\cref{prop:LR_convexity}) but not known to be $\triangle_{\min}$-star-quasi-convex, so does not meet the assumptions of \cref{prop:convergence_rate}; despite this, we see $\epsilon_m'$ is often an effective upper bound on the error in $\vec t$.}

\begin{figure}[h!]
	\centering
	\includegraphics[width=0.9\linewidth]{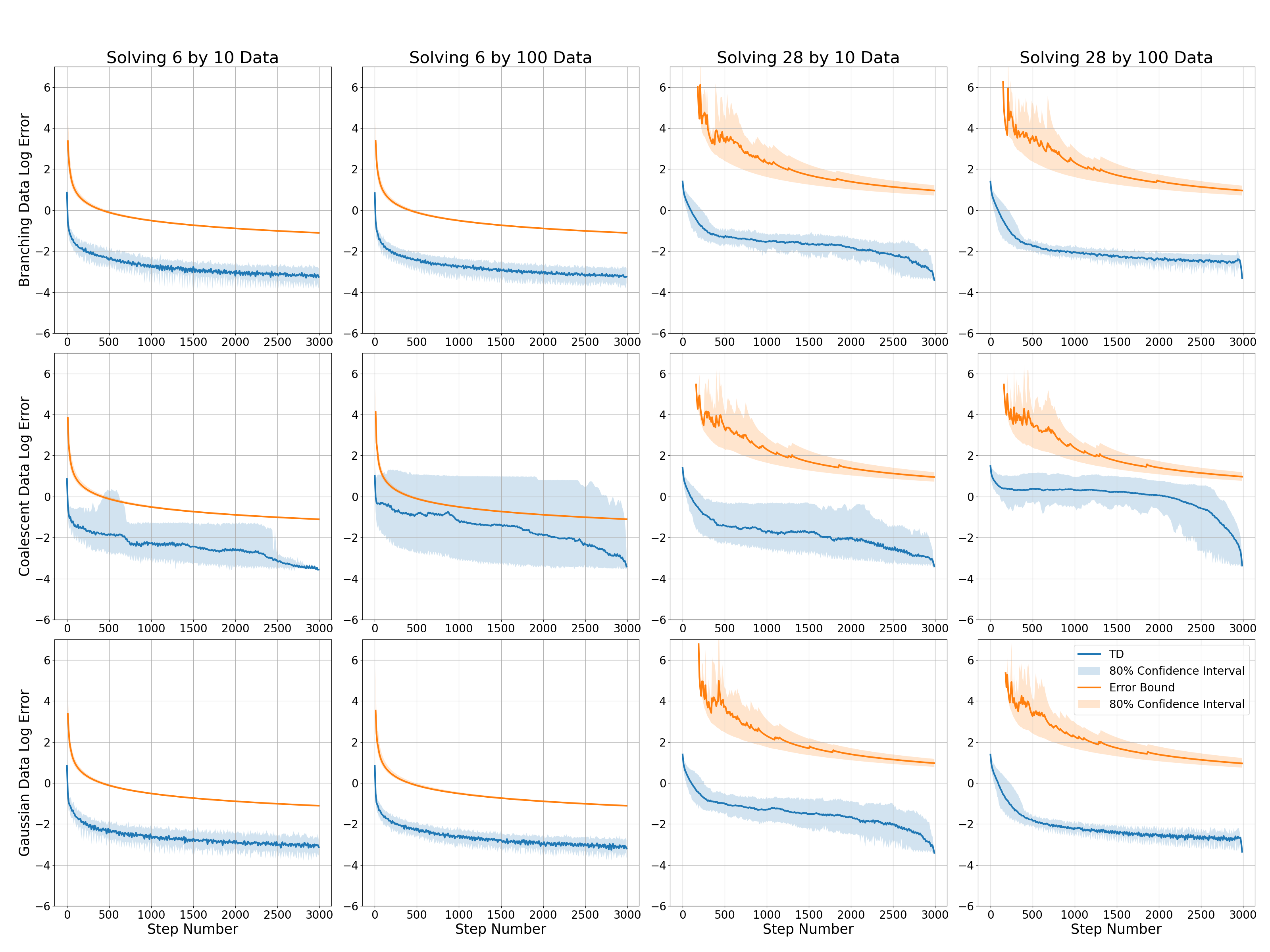}
	\caption{The mean log error bound $\log \, \epsilon_m'$ and the mean log tropical distance to the estimated minimum $\vec t^*$ when applying tropical descent (TD) with $\alpha = 1$ to solve the tropical linear regression problem. For each dataset, we take 50 random initialization points and compute $\log \, \epsilon_m'$ and $\log \, \dtr(\vec t_m, \vec t^*)$ independently, then take an average. We also record the 10th and 90th percentile values to estimate an 80\% confidence interval.}
	\label{fig:convergence_rate}
\end{figure}

{We perform tropical descent for tropical linear regression using 12 different datasets; these are generated via branching processes, coalescent processes, or a Gaussian distribution on $\TPT{6}$ and $\TPT{28}$. We consider sample sizes of 10 or 100. These are the same datasets used in \cref{sec:experiments}, and further details on their generation can be found in \cref{appsec:methodology}. We run computations using a learning rate of $\alpha = 1$, for steps up to $m \leq 3000$. The tropical linear regression problem has $\|\nabla f(\vec t_m)\| = 2$ almost everywhere, so we take $L=1$ when computing $\epsilon_m'$. We estimate $ \vec t^*$ to be the mean of the final 10 steps of tropical descent (TD), and compute a tight estimate for $D$ by taking 
	\[
	D\approx \max_{m \leq 3000} d_{\triangle_{\max}} (\vec t^*,\vec t_m).
	\]
	\cref{fig:convergence_rate} shows the mean log error bound $\epsilon_m'$ with the mean log tropical distance to $\vec t^*$.}

{Examining \cref{fig:convergence_rate}, we first note that the error $\epsilon_m'$ is orders of magnitude smaller for the lower dimensional problems, which is primarily a result of the constant terms in the denominator of $\epsilon_m'$. As we have noted, it would be of significant value to reduce these constant terms in future numerical analytical studies of tropical descent. These constant terms are also the cause of noise in the mean error bound for high dimensional problems, as we require more steps for the error to be well-defined. We also note that although $\epsilon_m'$ is a reasonable upper bound for the distance to $\vec t^*$ in many cases, this upper bound does not hold for the low dimensional coalescent datasets. It is unclear whither this is due to poor estimates for $D$ and $\vec t^*$, or due to tropical linear regression not being $\triangle_{\min}$-star-quasi-convex. We conclude by noting that we are observing a convergence rate of $O(1/\sqrt{m})$ for most experiments, with this rate appearing particularly distinctly for branching process and Gaussian data.}

\section{Numerical Experiments}
\label{sec:experiments}

In this section we compare the practical performance of classical descent, tropical descent, SGD, tropical SGD, Adam, Adamax, and {tropical Adamax}, with a focus on their stability behavior. The implementation for all these experiments is available in the GitHub repository at \url{https://github.com/Roroast/TropicalGradDescent}.

\subsection{Methods}

We begin by outlining our implementation of steepest descent, stochastic gradient descent (SGD), and Adam as well as their tropical variants. We take a heuristic approach to the tropicalization of SGD and Adam, replacing derivatives with tropical steepest descent directions in each case.

\textbf{Steepest Descent.}
We will perform classical descent (CD) and tropical descent (TD) using equivalent step size sequences $a_m$ for both methods, which we fix as
\[
a_{ m} = \frac{\gamma\| \nabla f(\vec t_{ m}) \|}{\sqrt{ m}},
\]
where $\| \cdot\|$ is the Euclidean norm for classical descent and the tropical norm for tropical descent. When $\nabla f(\vec t)$ is bounded this is a well-behaved step rule; $a_{ m}$ will converge to 0. Furthermore, for the piecewise linear problems, $\|\nabla f(\vec t_{ m})\|$ is bounded away from 0 almost everywhere which is not a global minimum, {so $\sum a_{m} = \infty$}.

\textbf{Stochastic Gradient Descent.}
Stochastic gradient descent is a variant of gradient descent which looks to improve computational speed by approximating $\nabla f(\vec t)$; rather than computing the loss and its derivative at $\vec t$ with respect to the entire dataset $f(\vec t; X)$, we evaluate the loss and its derivative with respect to a random sample $\vec x_{ m} \in X$, $f(\vec t; \vec x_{ m})$. This is simplified with the notation $f_{ m}(\vec t)$. As with gradient descent, we use the step size sequence given by
\[
a_{ m} = \frac{\gamma\| \nabla f_{ m}(\vec t_{ m}) \|}{\sqrt{ m}},
\]
where $\| \cdot\|$ is the Euclidean norm for stochastic gradient descent (SGD) and the tropical norm for tropical stochastic gradient descent (TSGD).

\textbf{Adam Variants.}
The Adam algorithm \citep{kingma2017adam} uses gradient momentum to compound consistent but small effects. This is done by taking first and second moments of the gradient over all steps using exponentially decaying weights. The Adamax algorithm is a variant of the Adam algorithm motivated by the use of the $\infty$-moments rather than second moments.

To adapt the Adam algorithm for the tropical setting, we take moments of an unnormalized tropical descent direction. When looking to take higher moments of the gradient over $ m$, the $L_{\infty}$ is the natural option for the tropical setting. We therefore look to use the Adamax recursion relation on the $\infty$-moment estimate $\vec u_{ m}$:
\[
\vec u_{ m+1} = \max(\beta_2 \vec u_{ m}, \vec d_{ m}).
\]

The tropicalized Adamax algorithm is then given by \Cref{alg:TrAdamax}.

\begin{algorithm}[ht!]
	\caption{TrAdamax}\label{alg:TrAdamax}
	\begin{algorithmic}
		\State ${\vec v}, \vec u \gets \vec 0$;
		\For{$ m < $ max\_steps}
		\State $\vec g_{ mi} \gets \nabla f(\vec t_{ m-1})_i$;
		\State $\vec d_{ mi} \gets (\max_j(\vec g_{ mj}) - \min_j(\vec g_{ mj})) * \mathbf{1}_{\vec g_{ mi} > 0}$;
		\State ${\vec v_{ m}} \gets \beta_1 * {\vec v_{ m-1}} + (1-\beta_1) * \vec d_{ m}$;
		\State $\vec u_{ m} \gets \max(\beta_2 * \vec u_{ m-1}, |\vec d_{ m}|)$;
		\State $\vec t_{ m} \gets \vec t_{ m} - \alpha/(1-\beta_1^{ m}) * {\vec v_{ m}}/\vec u_{ m}$;
		\EndFor
	\end{algorithmic}
\end{algorithm}

\subsection{Computation}

We now look to solve our centrality statistics, Wasserstein projection problems, and tropical linear regression using the 7 gradient methods discussed: classical descent (CD), tropical descent (TD), stochastic gradient descent (SGD), tropical stochastic gradient descent(TSGD), Adam, Adamax, and TrAdamax. We run computations for 12 datasets of different sizes, shapes and dimensions; motivated by phylogenetic data science applications, we consider datasets generated by branching processes, coalescent processes, or Gaussian distributions on $\TPT{N}$, for $N = 6, 28$ and samples of size 10 or 100. We take the same random initializations for each gradient method, taken from a Gaussian distribution on $\R^{N}$. \Cref{appsec:methodology} provides further details on the simulated datasets and computational methodology used. For each gradient method, we use learning rates which have been optimized for each dataset size and dimension, across different distribution types (\Cref{appsec:lr_tuning}). 

Rather than investigating the relative convergence rates of each method, we are interested in the relative stability of local minima. We therefore plot the CDF of the relative log error after 1000 steps over 50 random initializations, {to reveal concentrations of error values where multiple initializations stabilize}. When the CDF of one method dominates another, it has a higher probability of achieving the same error {or less} and is therefore a more suitable method.

As we will observe, tropical gradient methods are stable at fewer points than classical gradient methods; as a result, tropical methods typically have further to travel before reaching a stable point. Over few steps, this would be compensated for by selecting larger learning rates for tropical gradient methods relative to classical gradient methods, leading to less precise results in the case of tropical gradient methods. Our computations are run over 1000 steps, as this gives enough time for tropical gradient methods to {stabilize}, given our {initializations}, without sacrificing learning rate; we see in \Cref{appsec:lr_tuning} that neither classical nor tropical learning rates are consistently greater than the other.

\subsubsection{Centrality Statistics}\label{subsubsec:CS}

We first consider the behavior of tropical descent methods {for Fermat--Weber points and Fr\'echet means}. These optimization problems are classically convex problems, as well as being both $\max$-tropical and $\min$-tropical location problems, and therefore lend themselves well to gradient methods. The experimental behavior of these problems is similar, so in this section we include the results for the Fermat--Weber problem while the results for Fr\'echet means can be found in \Cref{appsec:further_figures}.

\begin{figure}[!pht]
	\centering
	\includegraphics[width = 0.95\textwidth]{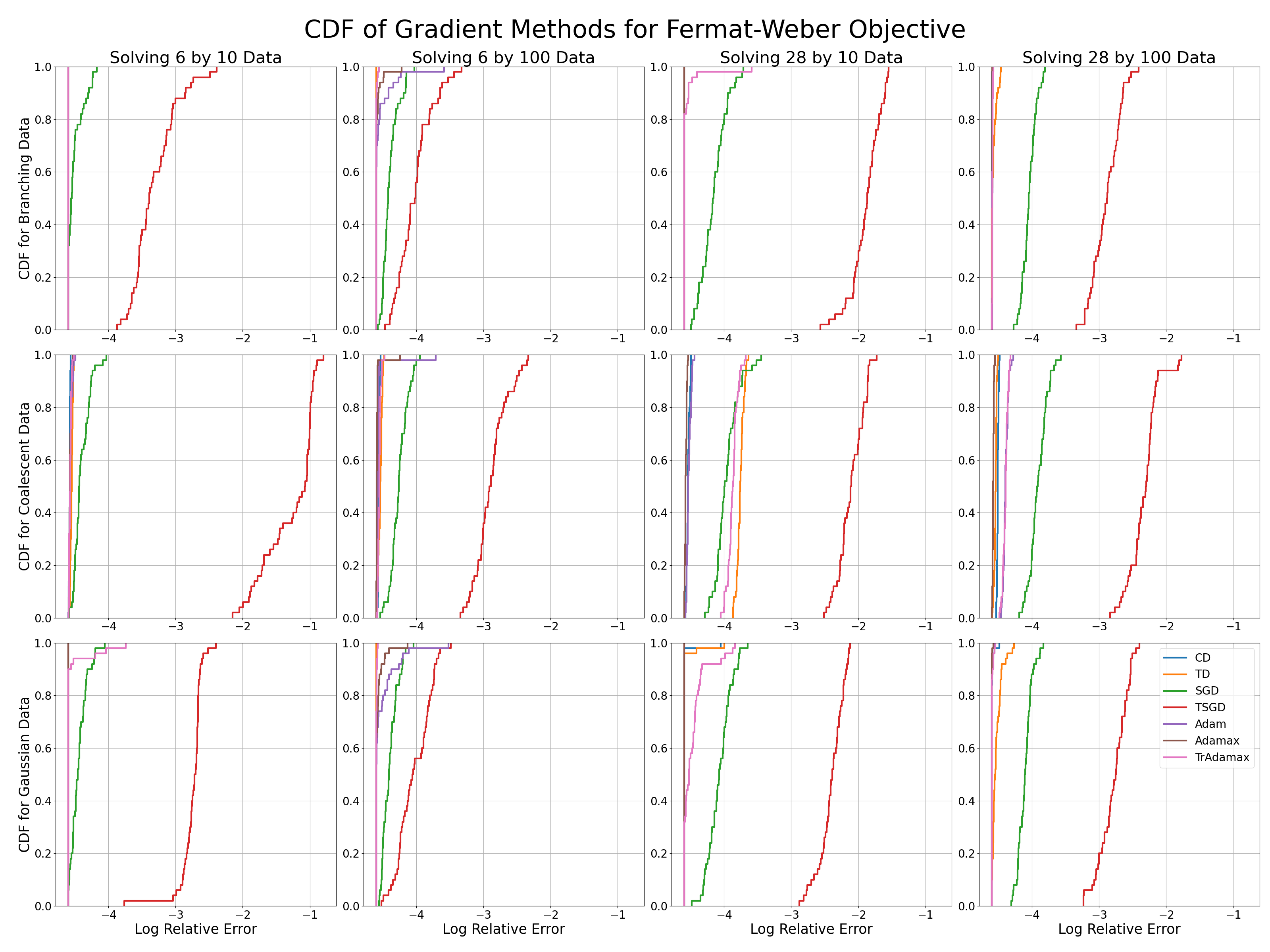}
	\caption{The CDF of the log relative error after 1000 steps for each gradient method with tuned learning rates across 50 random (Gaussian) initializations when minimizing the Fermat--Weber objective. Each subfigure corresponds to a dataset with sample size $K=10,100$ of dimensionality $N=6,28$, sampled from a branching process, coalescent process, or a Gaussian distribution.}
	\label{fig:FW_CDF_plot}
\end{figure}

\Cref{fig:FW_CDF_plot} shows the CDF of the log relative error {after} 1000 steps of each gradient method for the 12 datasets. We note that the deterministic methods achieve similar results for most datasets, while the stochastic methods achieve the worst error for all 12 datasets and TSGD is performing worse than SGD. The disparity between SGD and TSGD is greatest when the dimensionality of the data is greater than the size of the dataset, and for coalescent data. We see similar behavior in the log error of {tropical descent} and TrAdamax; they perform on par with their classical counterparts for each dataset other than the 28 by 10 coalescent data, for which the tropical gradient methods are not achieving the same accuracy (\cref{tab:FW_mean_log_errors}). While it is unclear why this is, we note that coalescent data is supported on the space of ultrametric trees, which is a tropical hyperplane \citep{ardila2006bergman,page2020tropical} and may influence the {stability} of minima. 

\begin{table}[!pht]
	\caption{The mean log relative error of each gradient method after 1000 steps with tuned learning rates for branching process, coalescent process, and Gaussian datasets of size $K=10,100$ in $\TPT{6}$, $\TPT{28}$ when minimizing the Fermat--Weber objective function. The minimal mean errors for each dataset are in bold.}
	\vspace{.1 in}
	\centering
	\begin{tabular}{|c|c|c|c|c|c|c|c|}
		\hline
		Data & CD & TD & SGD & TSGD & Adam & Adamax & TrAdamax \\
		\hline
		6$\times$10 Branching Data & \textbf{-4.60} & \textbf{-4.60} & -4.50 & -3.32 & \textbf{-4.60} & \textbf{-4.60} & \textbf{-4.60} \\
		6$\times$10 Coalescent Data & \textbf{-4.58} & -4.55 & -4.40 & -1.28 & -4.57 & -4.57 & -4.57 \\
		6$\times$10 Gaussian Data & \textbf{-4.60} & \textbf{-4.60} & -4.44 & -2.75 & \textbf{-4.60} & \textbf{-4.60} & -4.56 \\
		6$\times$100 Branching Data & \textbf{-4.60} & \textbf{-4.60} & -4.39 & -4.02 & -4.54 & -4.58 & -4.59 \\
		6$\times$100 Coalescent Data & -4.55 & -4.54 & -4.25 & -2.88 & -4.56 & \textbf{-4.58} & -4.55 \\
		6$\times$100 Gaussian Data & \textbf{-4.60} & \textbf{-4.60} & -4.39 & -4.03 & -4.52 & -4.57 & -4.59 \\
		28$\times$10 Branching Data & \textbf{-4.60} & \textbf{-4.60} & -4.15 & -1.89 & \textbf{-4.60} & \textbf{-4.60} & -4.56 \\
		28$\times$10 Coalescent Data & -4.53 & -3.76 & -3.96 & -2.11 & -4.52 & \textbf{-4.57} & -3.87 \\
		28$\times$10 Gaussian Data & -4.58 & -4.58 & -4.07 & -2.40 & \textbf{-4.60} & \textbf{-4.60} & -4.47 \\
		28$\times$100 Branching Data & \textbf{-4.59} & -4.57 & -4.04 & -2.89 & \textbf{-4.59} & \textbf{-4.59} & \textbf{-4.59} \\
		28$\times$100 Coalescent Data & -4.51 & -4.54 & -3.91 & -2.33 & -4.39 & \textbf{-4.58} & -4.39 \\
		28$\times$100 Gaussian Data & \textbf{-4.59} & -4.53 & -4.10 & -2.78 & \textbf{-4.59} & \textbf{-4.59} & \textbf{-4.59} \\
		\hline
	\end{tabular}
	\label{tab:FW_mean_log_errors}
\end{table}

\subsubsection{Wasserstein Projections}\label{subsubsec:WP}

\begin{figure}
	\centering
	\includegraphics[width = 0.95\textwidth]{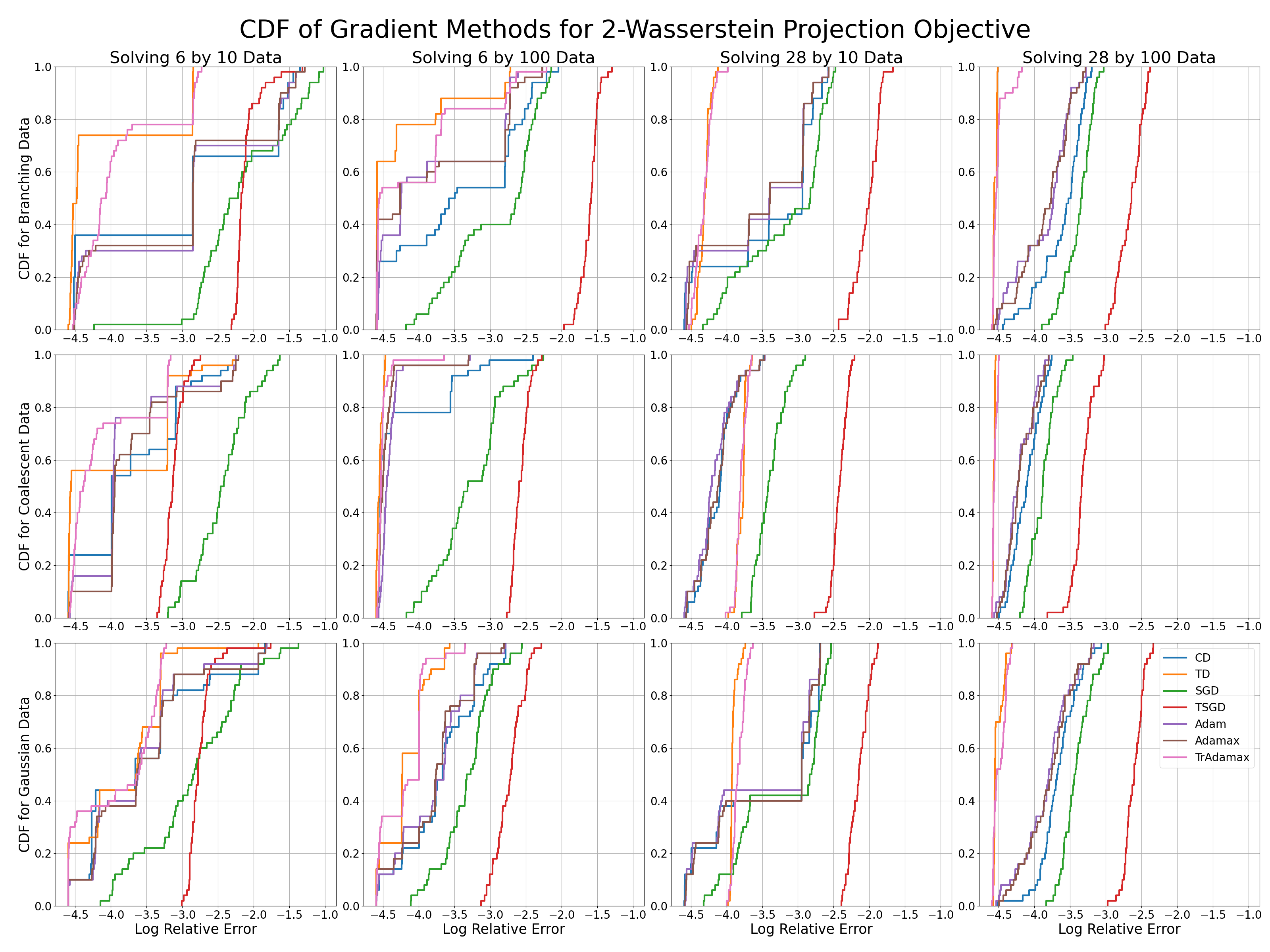}
	\caption{The CDF of the log relative error after 1000 steps for each gradient method with tuned learning rates across 50 random (Gaussian) initializations when minimizing the 2-Wasserstein projection objective. Each subfigure corresponds to a dataset with sample size $K=10,100$ of dimensionality $N=6,28$, sampled from a branching process, coalescent process, or a Gaussian distribution.}
	\label{fig:WD2_CDF_plot}
\end{figure}

\begin{figure}
	\centering
	\includegraphics[width = 0.95\textwidth]{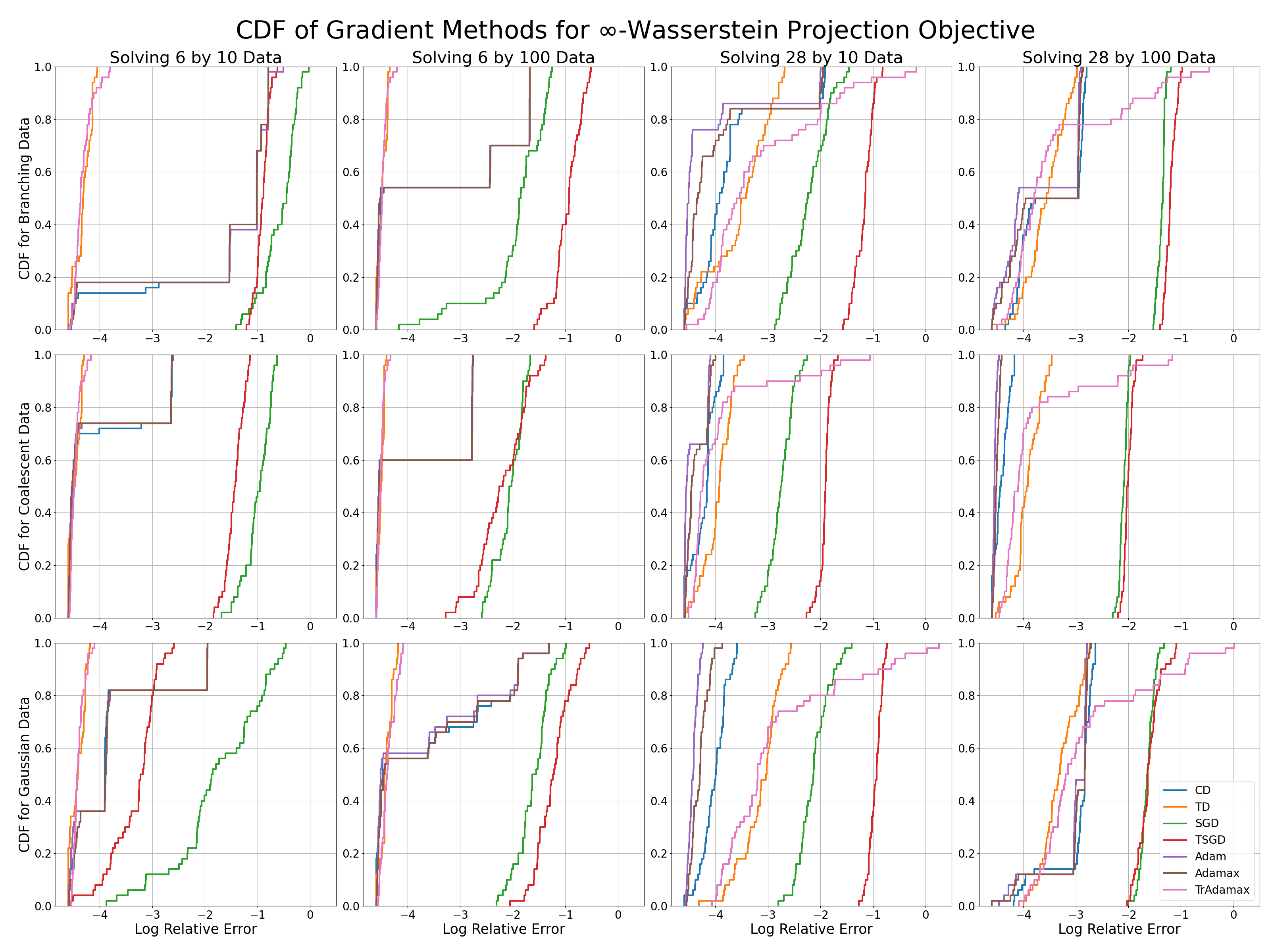}
	\caption{The CDF of the log relative error after 1000 steps for each gradient method with tuned learning rates across 50 random (Gaussian) initializations when minimizing the $\infty$-Wasserstein projection objective. Each subfigure corresponds to a dataset with sample size $K=10,100$ of dimensionality $N=6,28$, sampled from a branching process, coalescent process, or a Gaussian distribution.}
	\label{fig:WDinf_CDF_plot}
\end{figure}

The Wasserstein projection problems highlight the differences in {stability} behavior between classical and tropical gradient methods, as they are $\min$-tropical location problems but not classically convex. The CDFs of the final log relative errors can be found in \Cref{fig:WD2_CDF_plot,fig:WDinf_CDF_plot}.

In \Cref{fig:WD2_CDF_plot,fig:WDinf_CDF_plot} we see large steps at various intervals in the CDF plots, which is indicative of stable local minima. This is particularly prevalent in low dimensional space. Classical descent, Adam, and Adamax are particularly susceptible to this, while tropical descent and TrAdamax are more likely to pass these local minima. This difference is most pronounced when solving the $\infty$-Wasserstein problem in low dimensions---we are almost certain to reach the minimum with {relative error $\leq e^{-3}$}, while classical methods have a 50-80\% chance of achieving the same error.

\subsubsection{Linear Regression}\label{subsubsec:LR}

\begin{figure}[!pht]
	\centering
	\includegraphics[width = 0.95\textwidth]{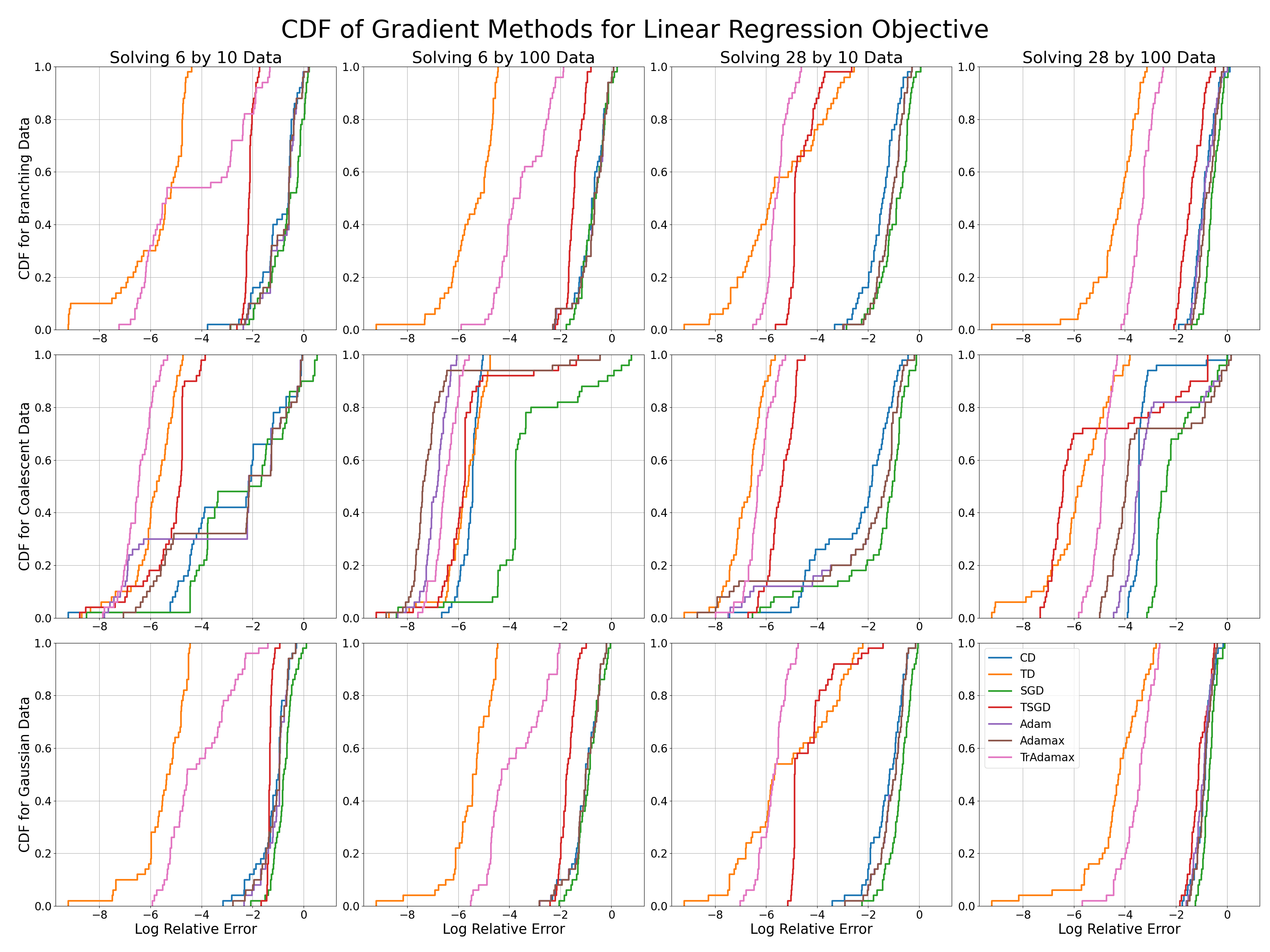}
	\caption{The CDF of the log error after 1000 steps for each gradient method with tuned learning rates across 50 random (Gaussian) initializations when minimizing the linear regression objective. Each subfigure corresponds to a dataset with sample size $K=10,100$ of dimensionality $N=6,28$, sampled from a branching process, coalescent process, or a Gaussian distribution.}
	\label{fig:LR_CDF_plot}
\end{figure}

As the coalescent data lies on a tropical hyperplane \citep{ardila2006bergman,page2020tropical}, the minimum of the linear regression function for such data is exactly 0. \Cref{fig:LR_CDF_plot} therefore shows the CDF of the log error rather than the CDF of the log relative error.

\begin{table}[!pht]
	\caption{The mean log relative error of each gradient method after 1000 steps with tuned learning rates for branching process, coalescent process, and Gaussian datasets of size $K=10,100$ in $\TPT{6}$, $\TPT{28}$ when minimizing the linear regression objective function. The minimal mean errors for each dataset are in bold.}
	\vspace{.1 in}
	\centering
	\begin{tabular}{|c|c|c|c|c|c|c|c|}
		\hline
		Data & CD & TD & SGD & TSGD & Adam & Adamax & TrAdamax \\
		\hline
		6$\times$10 Branching Data & -0.97 & \textbf{-5.79} & -0.56 & -2.16 & -0.80 & -0.84 & -4.39 \\
		6$\times$10 Coalescent Data & -2.70 & -5.87 & -2.11 & -5.47 & -2.97 & -2.66 & \textbf{-6.58} \\
		6$\times$10 Gaussian Data & -1.19 & \textbf{-5.48} & -0.80 & -1.35 & -1.05 & -1.10 & -4.10 \\
		6$\times$100 Branching Data & -0.79 & \textbf{-5.49} & -0.71 & -1.43 & -0.72 & -0.72 & -3.53 \\
		6$\times$100 Coalescent Data & -5.54 & -5.76 & -3.11 & -5.57 & -6.90 & \textbf{-7.01} & -6.52 \\
		6$\times$100 Gaussian Data & -1.04 & \textbf{-5.50} & -0.88 & -1.70 & -1.02 & -1.02 & -3.81 \\
		28$\times$10 Branching Data & -1.49 & -5.51 & -0.93 & -4.49 & -1.14 & -1.14 & \textbf{-5.57} \\
		28$\times$10 Coalescent Data & -2.43 & \textbf{-6.63} & -1.77 & -5.55 & -2.16 & -2.29 & -6.18 \\
		28$\times$10 Gaussian Data & -1.26 & -5.16 & -0.79 & -4.29 & -1.08 & -1.08 & \textbf{-5.72} \\
		28$\times$100 Branching Data & -0.90 & \textbf{-4.37} & -0.61 & -1.42 & -0.85 & -0.77 & -3.32 \\
		28$\times$100 Coalescent Data & -3.33 & \textbf{-5.84} & -2.17 & -5.35 & -3.02 & -3.20 & -4.93 \\
		28$\times$100 Gaussian Data & -0.91 & \textbf{-4.35} & -0.62 & -1.07 & -0.96 & -0.90 & -3.51 \\
		\hline
	\end{tabular}
	\label{tab:LR_mean_log_errors}
\end{table}

\Cref{tab:LR_mean_log_errors} shows the mean log errors of each gradient method and dataset, from which we see that the tropical methods outperform their classical counterparts for all but one dataset, often by a factor of two. 
In \Cref{fig:LR_CDF_plot} we see that the CDF of log errors for tropical descent and TrAdamax consistently dominates those of classical descent, Adam and Adamax, to the extent that the support of their CDFs does not intersect that of their classical counterparts for some datasets.
We also note that unlike the previous problems, TSGD is outperforming SGD when solving the linear regression problem, often to the extent that TSGD outperforms all classical gradient methods.

{\subsection{Tropical Statistics for Applications}}

{We conclude our computational experiments by using tropical gradient methods to solve tropical statistical problems which appear in phylogenetic data analysis and game theory. In particular we look to use tropical centrality statistics to estimate the species tree of a multi-species coalescent model \citep{liu2019modern} , and we use tropical linear regression to estimate the hidden preference factors in the Auction model introduced by \cite{akian2023tropical}, which is a variant of the first-price sealed-bid auction model \citep{krishna2009auction}.}

{\subsubsection{Multi-species Coalescent Model}}

{A \emph{species tree} refers to a rooted metric tree which represents the evolutionary relationships between a set of species. However, for a specific gene, we may not observe the same branching pattern, instead observing a gene-specific evolution which is referred to as a \emph{gene tree}. This raises one of the key questions in phylogenetics; how do we estimate the species tree from a set of gene trees? The multi-species coalescent model is a probabilistic model for gene trees given a certain species tree (see \cite{liu2019modern} for a detailed introduction), which is particularly valuable for formulating Bayesian estimators for species trees and producing simulated data with which to test species tree estimation methodologies.}

\begin{figure}
	\centering
	\includegraphics[width=0.5\linewidth]{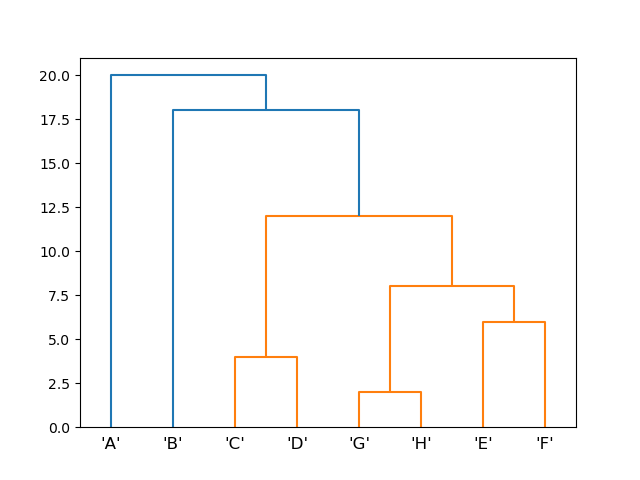}
	\caption{The underlying species tree used for the multi-species coalescent model generating our dataset of gene trees. The topology of this tree can be represented by the Newick string ``(((((F,E),(H,G)),(D,C)),B),A)''.}
	\label{fig:sp_tree}
\end{figure}

{In this section, we consider the tropical Fermat-Weber and Fr\'echet means of gene trees as estimators for the underlying species tree. We fix a species tree with 8 leaves (see \cref{fig:sp_tree}), then use Dendropy \citep{dendropy5} to sample 100 gene trees from a multi-species coalescent model with this underlying species tree and normalize this dataset to have average tropical norm 1.}

{We then compute a Fermat-Weber point and Fr\'echet mean for this dataset using each of our gradient methods for 1000 steps with their trained learning rates (\cref{appsec:lr_tuning}). We perform these computations from 100 random initializations for each gradient method, then project the result to ultrametric space using the single-linkage algorithm, which is known to coincide with tropical projection \citep{single-linkage,ardila2005subdominant}, and record the topology of the resulting tree. The counts for each observed tree topology are presented in \cref{tab:MSC_FW_counts,tab:MSC_FM_counts}.}
\begin{table}[]
	\centering
	\begin{tabular}{|c|c|c|c|c|c|c|c|}
		\hline
		Newick Tree & Adam & Adamax & CD & SGD & TD & TrAdamax & TSGD \\
		\hline
		((B,A),(((F,E),(H,G)),(D,C))); & 0 & 0 & 0 & 0 & 0 & 0 & 2 \\
		(((((F,E),(H,G)),(D,C)),A),B); & 0 & 0 & 0 & 3 & 0 & 0 & 0 \\
		(((((F,E),(H,G)),(D,C)),B),A); & 100 & 100 & 100 & 97 & 100 & 100 & 98 \\
		\hline
	\end{tabular}
	\caption{The counts of each Newick string (up to tree topology) returned by the tropical projection of a Fermat--Weber point computed via each gradient method. Each method is run for 1000 steps with trained learning rate, from 100 random initializations.}
	\label{tab:MSC_FW_counts}
\end{table}
\begin{table}[]
	\centering
	\begin{tabular}{|c|c|c|c|c|c|c|c|}
		\hline
		Newick Tree & Adam & Adamax & CD & SGD & TD & TrAdamax & TSGD \\
		\hline
		((B,A),(((F,E),(H,G)),(D,C))); & 0 & 0 & 0 & 0 & 0 & 0 & 5 \\
		(((((F,E),(H,G)),(D,C)),A),B); & 0 & 0 & 0 & 1 & 0 & 0 & 1 \\
		(((((F,E),(H,G)),(D,C)),B),A); & 100 & 100 & 100 & 99 & 100 & 100 & 94 \\
		\hline
	\end{tabular}
	\caption{The counts of each Newick string (up to tree topology) returned by the tropical projection of a Fr\'echet means computed via each gradient method. Each method is run for 1000 steps with trained learning rate, from 100 random initializations.}
	\label{tab:MSC_FM_counts}
\end{table}

{We first note that in the vast majority of cases, our calculations recover the tree topology of the original species tree, suggesting that for this multi-species coalescent dataset, the tropical Fermat--Weber and Fr\'echt mean points are sufficiently good estimators for the species tree. When comparing the performance of each gradient method, we notice that the stochastic methods are the only methods to not recover the species tree topology, with Fr\'echet means computed via TSGD producing the highest error rate.}

{The relative performance of tropical Fermat--Weber and Fr\'echet mean points as estimators for multi-species coalescent species trees is currently an open topic of study \citep{cox2025tropical}. In exploring the consistency of these estimators it world be natural for future work to also consider the relative performance of tropical gradient methods for their computation.}

{\subsubsection{Auction Model with Hidden Preference Factors}}\label{subsubsec:Auction}

{As a final computational demonstration, we use tropical gradient methods to solve the problem of tropical linear regression using the toy auction data example given by \cite{akian2023tropical}. For $N$ firms and $K$ products, we let $(p_{ij})_{i \leq N, j \leq K}$ denote the price firm $i$ offers for product $j$. Assuming the offered prices are in a state of equilibrium, that is, no firm can raise their price and still win, differences in proposed price from one firm relative to another is indicative of a hidden preference factor $f_i$ in the market; that is, two firms' offers for product $j$ are considered equal when $p_{i_1j}f_{i_1} = p_{i_2j}f_{i_2}$. By performing tropical linear regression on the data $V_{ij} = \exp(-p_{ij})$, we can detect these hidden preference factors. For further details of this model and its assumptions, see \cite{akian2023tropical}.}

{For our computations, we take $N=3$ and $K=6,100$, using the same toy datasets discussed in Section 7.2 of \cite{akian2023tropical}, which have true hidden preference factors $f = (1, 0.8, 0.6)$. We then perform each gradient method for 1000 steps over 100 random initializations, using learning rates which have been tuned to these datasets (\cref{appsec:lr_tuning}). Our estimates for $f$ are then given by $\hat f \coloneqq \exp(\vec t_{1000}-\max(\vec t_{1000}))$, such that $\max_i \hat f = 1$.}

\begin{table}[]
	\centering
	\begin{tabular}{|c|cc|}
		\hline
		& K=6 & K=100 \\
		\hline
		Adam & (0.76 $\pm$ 0.33, 0.70 $\pm$ 0.32, 0.56 $\pm$ 0.21) & (0.77 $\pm$ 0.33, 0.71 $\pm$ 0.32, 0.58 $\pm$ 0.21) \\
		Adamax & (0.76 $\pm$ 0.33, 0.70 $\pm$ 0.32, 0.56 $\pm$ 0.21) & (0.77 $\pm$ 0.33, 0.71 $\pm$ 0.32, 0.59 $\pm$ 0.21) \\
		CD & (0.71 $\pm$ 0.43, 0.60 $\pm$ 0.42, 0.45 $\pm$ 0.30) & (0.78 $\pm$ 0.33, 0.70 $\pm$ 0.32, 0.58 $\pm$ 0.20) \\
		SGD & (0.73 $\pm$ 0.44, 0.57 $\pm$ 0.43, 0.41 $\pm$ 0.32) & (0.77 $\pm$ 0.33, 0.70 $\pm$ 0.32, 0.58 $\pm$ 0.21) \\
		TD & (1.00 $\pm$ 0.00, 0.81 $\pm$ 0.02, 0.60 $\pm$ 0.01) & (1.00 $\pm$ 0.02, 0.75 $\pm$ 0.08, 0.62 $\pm$ 0.05) \\
		TrAdamax & (1.00 $\pm$ 0.00, 0.81 $\pm$ 0.00, 0.60 $\pm$ 0.00) & (0.93 $\pm$ 0.19, 0.79 $\pm$ 0.16, 0.65 $\pm$ 0.08) \\
		TSGD & (1.00 $\pm$ 0.00, 0.82 $\pm$ 0.02, 0.60 $\pm$ 0.01) & (1.00 $\pm$ 0.04, 0.79 $\pm$ 0.04, 0.61 $\pm$ 0.01) \\
		\hline
	\end{tabular}
	\caption{The mean($\pm$ 1 std.) estimators $\hat f \coloneqq \exp(\vec t_{1000}-\max(\vec t_{1000}))$ for the hidden preference factors in the Auction model \citep{akian2023tropical} when computed using gradient methods to solve the tropical linear regression problem. The true hidden preference factors are $f= (1,0.8,0.6)$. \cite{akian2023tropical} uses the projective Krasnoselskii-Mann iteration to solve this tropical linear regression problem, which produced an estimate of $f^{\text{reg}} = (1, 0.81, 0.605)$ for the $K=6$ toy dataset.}
	\label{tab:Auction_Data}
\end{table}

{\cref{tab:Auction_Data} shows the resulting mean estimates for these hidden preference factors when using different gradient methods to solve the tropical linear regression problem. We see that the accuracy of the tropical gradient methods vastly exceeds each of the classical gradient methods, and that these gradient methods produce particularly consistent results for the small sample problem.}

{These computations are limited to a single Auction model, but show promising results. For future work, it would be of interest to explore how the precision of these methods changes with sample size and dimension, as well as a running time comparison with existing methods such as the projective Krasnoselskii-Mann iteration used by \cite{akian2023tropical} or combinatorial simplex algorithms \citep{allamigeon2014combinatorial}. \\
}

\section{Conclusion and Discussion}
\label{sec:conclusion}

We have introduced the concept of tropical descent as the steepest descent method with respect to the tropical norm. We have demonstrated, both theoretically and experimentally, that it {exhibits desirable} properties with respect to a wide class of optimization problems on the tropical projective torus, including several key statistical problems for the analysis of phylogenetic data. As well as proposing our own framework for tropical gradient methods, this work also provides a detailed review of the varying behavior of popular gradient methods when applied to tropical statistical optimization problems.

Our theoretical results outline the relative strengths of tropical descent over classical descent for tropical location problems, but only provide global solvability for 1-sample problems. In practice, our experimental results go {much} further; the optimization problems studied are location problems, rather than {$\triangle_{\min}$}-star-quasi-convex problems, and yet tropical descent appears to consistently converge to a global minimum, not just the tropical convex hull of our dataset. It would therefore be a natural extension of this work to identify how much we can relax the {$\triangle_{\min}$}-star-quasi-convexity condition while guaranteeing {global solvability}.

In this paper we have considered tree data with 4 to 8 leaves, and samples of size 10 to 100---at this scale, we have observed (T/)SGD to be a significant sacrifice in accuracy with little benefit in computational time. Furthermore, at this scale, the relative merits of SGD versus its tropical counterpart are inconsistent; TSGD out-performs SGD when solving the linear regression problem, but performs worse when finding centrality statistics. In a similar vein, initial investigations suggest that the classical stochastic algorithm for the computation of Fr\'echet means---Sturm's algorithm \citep{sturm2003probability}---performs poorly for the computation of tropical Fr\'echet means due to the inconsistent curvature of the tropical projective torus \citep{amendola2021invitation}. However, for applications to larger datasets and more complex objective functions, stochastic gradient methods may become necessary. It would then be important to understand the theoretical differences between SGD and TSGD in terms of their convergence guarantees in the tropical setting and their comparative performance at different scales of data dimension.

The recent identification of tropical rational functions and ReLU neural networks \citep{zhang2018tropical} has inspired new avenues of research in tropical applications, such as the tropical geometry of neural networks \citep{brandenburg2024real,yoshida2023tropNNs}. It is therefore natural to consider the possible extensions of this work to neural network training; can we identify local tropical quasi-convexity in the loss functions of ReLU neural networks, and how can we extend tropical descent to the space of tropical rational functions?

Finally, as shown by \cite{akian2023tropical}, the tropical linear regression problem is polynomial-time equivalent to mean payoff games---two-player perfect information games on a directed graph \citep{ehrenfeucht1979positional,gurvich1988cyclic}. This family of problems is of particular interest in game theory as they are known to be NP $\cap$ co-NP. In our numerical experiments, tropical gradient methods demonstrate a particularly strong performance in solving tropical linear regression, and while our methods are designed to approximate the minimum rather than compute it exactly, tropical descent may prove to be an efficient tool when solving large-scale mean payoff games.

\paragraph{Acknowledgements}
R.T.~receives partial funding from a Technical University of Munich--Imperial College London Joint Academy of Doctoral Studies (JADS) award (2021 cohort, PIs Drton/Monod).  A.M.~is supported by EPSRC AI Hub on Mathematical Foundations of Intelligence:
An ``Erlangen Programme'' for AI No. EP/Y028872/1.

\section*{Declarations}

The authors have no relevant financial or non-financial interests to disclose.
The authors have no competing interests to declare that are relevant to the content of this article.
All authors certify that they have no affiliations with or involvement in any organization or entity with any financial interest or non-financial interest in the subject matter or materials discussed in this manuscript.
The authors have no financial or proprietary interests in any material discussed in this article.

\appendix
	
\crefalias{section}{appendix}

\section{PCA and Logistic Regression} \label{appsec:further_experiments}

\begin{figure}[!hpt]
	\centering
	\includegraphics[width = 0.95\textwidth]{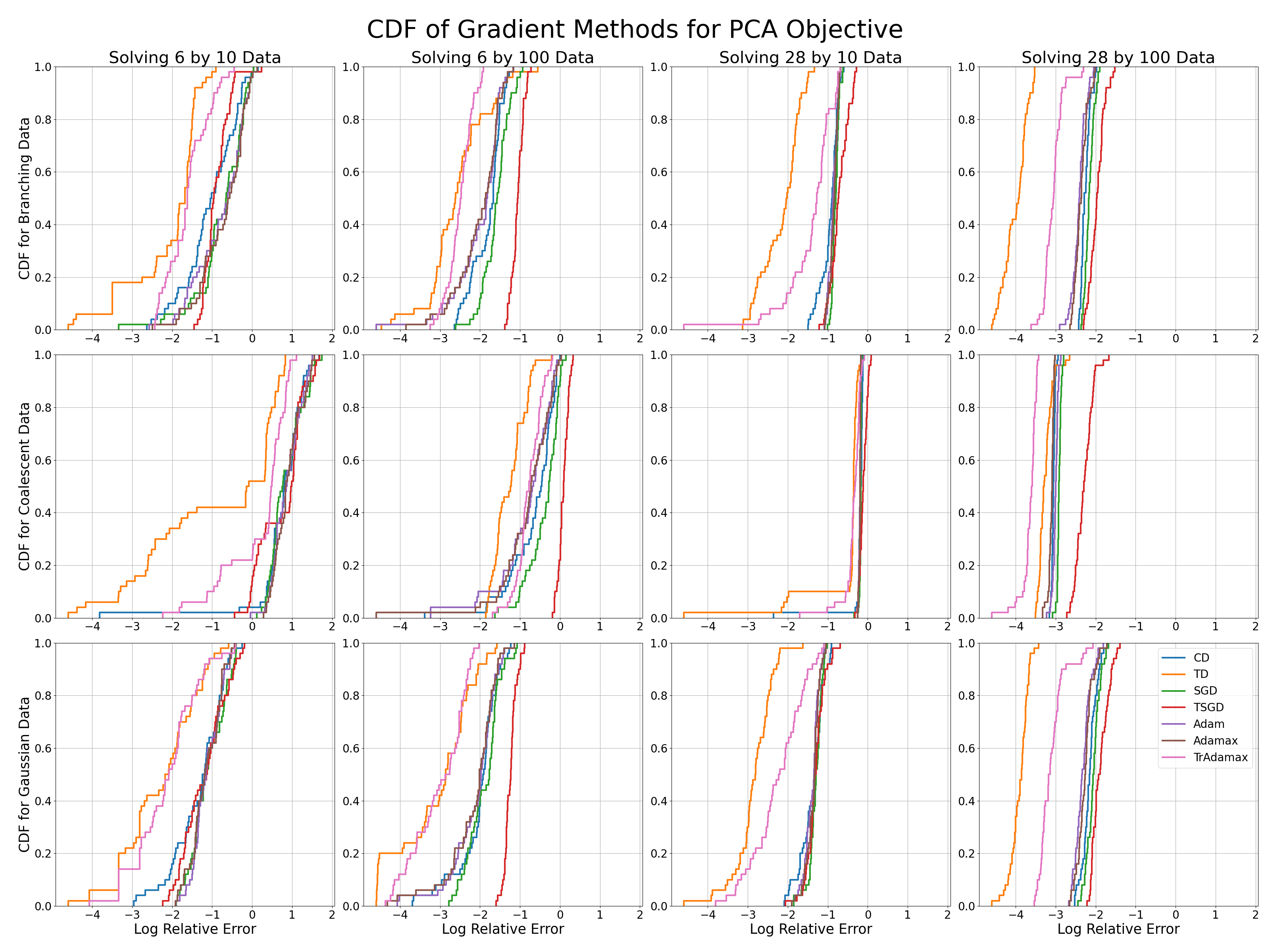}
	\caption{The CDF of the log relative error after 1000 steps for each gradient method with tuned learning rates across 50 random (Gaussian) initializations when minimizing the PCA objective. Each subfigure corresponds to a dataset with sample size $K=10,100$ of dimensionality $N=6,28$, sampled from a branching process, coalescent process, or a Gaussian distribution.}
	\label{fig:PCA_CDF_plot}
\end{figure}

\begin{figure}[!pht]
	\centering
	\includegraphics[width = 0.95\textwidth]{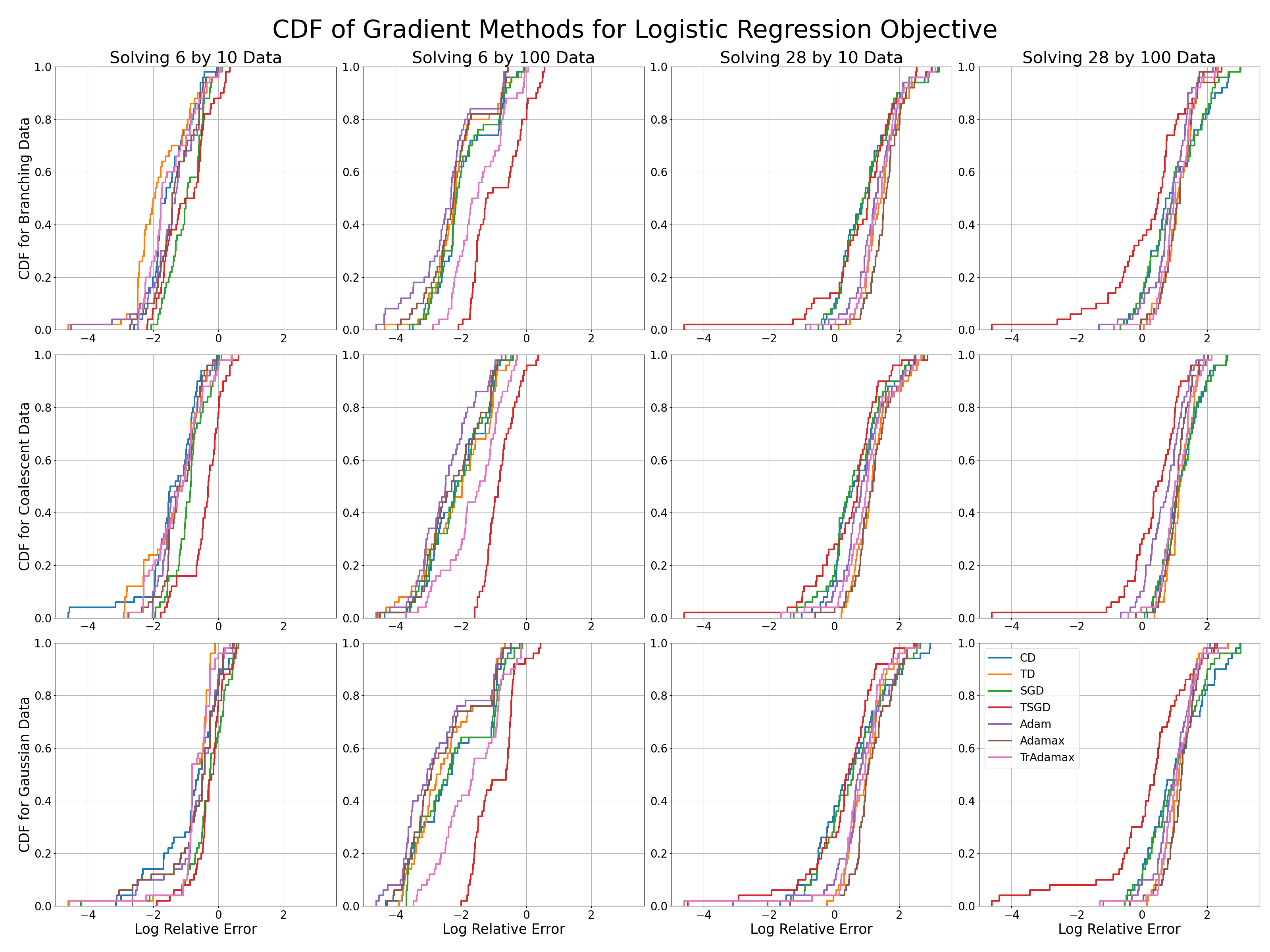}
	\caption{The CDF of the log relative error after 1000 steps for each gradient method with tuned learning rates across 50 random (Gaussian) initializations when minimizing the logistic regression objective. Each subfigure corresponds to a dataset with sample size $K=10,100$ of dimensionality $N=6,28$, sampled from a branching process, coalescent process, or a Gaussian distribution.}
	\label{fig:LogR_CDF_plot}
\end{figure}

As well as the location problems outlined in this paper, we performed the same numerical experiments with the tropical polytope PCA loss function \citep{yoshida2019tropicalPCA} and tropical logistic regression loss function \citep{aliatimis2023logistic}. These loss functions take multiple points in $\TPT{N}$ as argument, and we execute our gradient methods on each argument simultaneously. \Cref{fig:PCA_CDF_plot,fig:LogR_CDF_plot} show the resulting CDF of our gradient methods for these problems. These problems do not lend themselves well to gradient methods, as the relative log error is rather large for all of our proposed methods. However, for both problems the tropical descent and TrAdamax methods are outperforming the classical methods quite consistently.

\section{Methodology} \label{appsec:methodology}

The same 12 datasets are used for each experiment in \Cref{subsec:convergence,sec:experiments}. These datasets live in $\TPT{6}$ or $\TPT{28}$; these dimensions correspond to the spaces of trees with 4 and 8 leaves respectively. Branching process datasets are generated using the \texttt{R} package \texttt{ape} \citep{paradis2019ape} with exponentially distributed branch lengths---this produces data supported on the space of metric trees, i.e., the tropical Grassmannian. Coalescent data has also been generated using \texttt{ape} and is supported on the space of ultrametric trees, which is a tropical linear space. Gaussian data is generated using \texttt{NumPy} \citep{harris2020array}, and is supported on the general tropical projective torus.

In the computation of the Wasserstein projection objectives, we require a secondary dataset of different dimensionality; for each of the 12 datasets, we generate datasets of the same size and distribution, but in $\TPT{3}$ or $\TPT{21}$ which contain the spaces of trees with 3 and 7 leaves respectively.

We use \texttt{PyTorch} \citep{paszke2019pytorch} for the computation and differentiation of each objective function. Max-tropical descent directions are used for tropical gradient methods; we therefore optimize the min-tropical location problems---{Wasserstein} projections and tropical linear regression---with respect to $-\vec t$.

\section{Learning Rate Tuning}\label{appsec:lr_tuning}

The Adam algorithms---Adam, Adamax and TrAdamax---are run using the default $\beta = (0.9, 0.999)$ parameters. We then tune the learning rates used for each method; CD, TD, SGD, TSGD, Adam, Adamax, and TrAdamax. {We take 10 initializations at unit intervals on a log scale in the region $[e^{-6}, e^4]$ for each problem, gradient method, and data dimensionality. We then selected the learning rate with minimal mean log relative error in each case; \Cref{tab:FW_learning_rates,tab:FM_learning_rates,tab:WD2_learning_rates,tab:WDinf_learning_rates,tab:PCA_learning_rates,tab:LogR_learning_rates,tab:LR_learning_rates} show the selected learning rates.}

{For the Auction data in \cref{subsubsec:Auction}, we take 10 initializations at unit intervals in the region of $[e^{-10},e^5]$ for each dataset; \cref{tab:Auction_LR_learning_rates} shows the selected learning rates which minimize the mean log relative error.
}

\begin{table}[!pht]
	\caption{Tuned learning rates for each gradient method when solving the Fermat--Weber problem.}
	\vspace{.1 in}
	\centering
	\begin{tabular}{|l|l|l|l|l|l|l|l|l|}
		\hline
		\makecell{Data\\ Dim} & \makecell{Data\\ Count} & CD & TD & SGD & TSGD & Adam & Adamax & TrAdamax \\
		\hline
		\multirow[t]{2}{*}{6} & 10 & 0.135 & 0.135 & 0.135 & 0.135 & 0.00674 & 0.00674 & 0.00674 \\
		& 100 & 0.135 & 0.135 & 0.135 & 0.135 & 0.00248 & 0.00248 & 0.00674 \\
		\hline
		\multirow[t]{2}{*}{28} & 10 & 0.368 & 0.368 & 0.368 & 1 & 0.00674 & 0.00674 & 0.0183 \\
		& 100 & 0.368 & 0.135 & 0.368 & 1 & 0.00674 & 0.00674 & 0.0183 \\
		\hline
	\end{tabular}
	\label{tab:FW_learning_rates}
\end{table}

\begin{table}[!pht]
	\caption{Tuned learning rates for each gradient method when solving the Fr\'echet mean problem.}
	\vspace{.1 in}
	\centering
	\begin{tabular}{|l|l|l|l|l|l|l|l|l|}
		\hline
		\makecell{Data\\ Dim} & \makecell{Data\\ Count} & CD & TD & SGD & TSGD & Adam & Adamax & TrAdamax \\
		\hline
		\multirow[t]{2}{*}{6} & 10 & 0.135 & 0.135 & 0.135 & 0.135 & 0.00674 & 0.00674 & 0.00674 \\
		& 100 & 0.135 & 0.135 & 0.135 & 0.135 & 0.00674 & 0.00248 & 0.00674 \\
		\hline
		\multirow[t]{2}{*}{28} & 10 & 0.368 & 0.368 & 0.368 & 1 & 0.00674 & 0.00674 & 0.0183 \\
		& 100 & 0.368 & 0.135 & 0.368 & 1 & 0.00674 & 0.00674 & 0.0183 \\
		\hline
	\end{tabular}
	\label{tab:FM_learning_rates}
\end{table}

\begin{table}[!pht]
	\caption{Tuned learning rates for each gradient method when solving the 2-Wasserstein projection problem.}
	\vspace{.1 in}
	\centering
	\begin{tabular}{|l|l|l|l|l|l|l|l|l|}
		\hline
		\makecell{Data\\ Dim} & \makecell{Data\\ Count} & CD & TD & SGD & TSGD & Adam & Adamax & TrAdamax \\
		\hline
		\multirow[t]{2}{*}{6} & 10 & 2.72 & 1 & 1 & 0.135 & 0.135 & 0.368 & 0.368 \\
		& 100 & 2.72 & 2.72 & 0.368 & 0.135 & 0.135 & 0.135 & 0.135 \\
		\hline
		\multirow[t]{2}{*}{28} & 10 & 0.368 & 0.368 & 0.368 & 1 & 0.00674 & 0.0183 & 0.0183 \\
		& 100 & 2.72 & 2.72 & 0.368 & 1 & 0.0498 & 0.135 & 0.0498 \\
		\hline
	\end{tabular}
	\label{tab:WD2_learning_rates}
\end{table}

\begin{table}[!pht]
	\caption{Tuned learning rates for each gradient method when solving the $\infty$-Wasserstein projection problem.}
	\vspace{.1 in}
	\centering
	\begin{tabular}{|l|l|l|l|l|l|l|l|l|}
		\hline
		\makecell{Data\\ Dim} & \makecell{Data\\ Count} & CD & TD & SGD & TSGD & Adam & Adamax & TrAdamax \\
		\hline
		\multirow[t]{2}{*}{6} & 10 & 0.0498 & 0.135 & 0.0498 & 0.368 & 0.00674 & 0.00674 & 0.0183 \\
		& 100 & 0.0498 & 0.135 & 0.135 & 2.72 & 0.00674 & 0.00674 & 0.0183 \\
		\hline
		\multirow[t]{2}{*}{28} & 10 & 0.368 & 1 & 0.368 & 1 & 0.00674 & 0.0183 & 0.0498 \\
		& 100 & 0.368 & 1 & 0.368 & 1 & 0.00674 & 0.0183 & 0.0498 \\
		\hline
	\end{tabular}
	\label{tab:WDinf_learning_rates}
\end{table}

\begin{table}[!pht]
	\caption{Tuned learning rates for each gradient method when solving the polytope PCA problem.}
	\vspace{.1 in}
	\centering
	\begin{tabular}{|l|l|l|l|l|l|l|l|l|}
		\hline
		\makecell{Data\\ Dim} & \makecell{Data\\ Count} & CD & TD & SGD & TSGD & Adam & Adamax & TrAdamax \\
		\hline
		\multirow[t]{2}{*}{6} & 10 & 2.72 & 2.72 & 1 & 1 & 0.135 & 0.135 & 0.368 \\
		& 100 & 2.72 & 7.39 & 1 & 1 & 0.135 & 0.135 & 1 \\
		\hline
		\multirow[t]{2}{*}{28} & 10 & 7.39 & 7.39 & 0.368 & 1 & 0.0183 & 0.0183 & 0.135 \\
		& 100 & 2.72 & 7.39 & 0.368 & 1 & 0.0498 & 0.0498 & 0.0498 \\
		\hline
	\end{tabular}
	\label{tab:PCA_learning_rates}
\end{table}

\begin{table}[!pht]
	\caption{Tuned learning rates for each gradient method when solving the logistic regression problem.}
	\vspace{.1 in}
	\centering
	\begin{tabular}{|l|l|l|l|l|l|l|l|l|}
		\hline
		\makecell{Data\\ Dim} & \makecell{Data\\ Count} & CD & TD & SGD & TSGD & Adam & Adamax & TrAdamax \\
		\hline
		\multirow[t]{2}{*}{6} & 10 & 7.39 & 2.72 & 1 & 2.72 & 0.368 & 1 & 0.0498 \\
		& 100 & 1 & 0.368 & 1 & 2.72 & 0.00248 & 0.00248 & 0.00248 \\
		\hline
		\multirow[t]{2}{*}{28} & 10 & 1 & 0.368 & 1 & 2.72 & 0.00248 & 0.00248 & 0.00674 \\
		& 100 & 1 & 0.368 & 1 & 7.39 & 0.00248 & 0.00248 & 0.00674 \\
		\hline
	\end{tabular}
	\label{tab:LogR_learning_rates}
\end{table}

\begin{table}[!pht]
	\caption{Tuned learning rates for each gradient method when solving the linear regression problem.}
	\vspace{.1 in}
	\centering
	\begin{tabular}{|l|l|l|l|l|l|l|l|l|}
		\hline
		\makecell{Data\\ Dim} & \makecell{Data\\ Count} & CD & TD & SGD & TSGD & Adam & Adamax & TrAdamax \\
		\hline
		\multirow[t]{2}{*}{6} & 10 & 0.368 & 0.135 & 0.368 & 0.135 & 0.00248 & 0.0183 & 0.00674 \\
		& 100 & 0.135 & 0.135 & 0.368 & 0.0498 & 0.00248 & 0.00248 & 0.00674 \\
		\hline
		\multirow[t]{2}{*}{28} & 10 & 0.368 & 0.0498 & 0.0498 & 0.135 & 0.00248 & 0.00248 & 0.0183 \\
		& 100 & 1 & 0.368 & 1 & 7.39 & 0.0498 & 0.0498 & 0.0498 \\
		\hline
	\end{tabular}
	\label{tab:LR_learning_rates}
\end{table}

\begin{table}[!pht]
	\caption{Tuned learning rates for each gradient method when solving the linear regression problem with Auction data.}
	\vspace{.1 in}
	\centering
	\begin{tabular}{|l|l|l|l|l|l|l|l|l|}
		\hline
		\makecell{Data\\ Dim} & \makecell{Data\\ Count} & CD & TD & SGD & TSGD & Adam & Adamax & TrAdamax \\
		\hline
		\multirow[t]{2}{*}{3} & 6 & 2.72 & 0.368 & 7.39 & 0.0498 & 0.000912 & 0.000912 & 0.0183 \\
		& 100 & 0.0183 & 2.72 & 0.0183 & 0.0498 & 0.000912 & 0.000912 & 0.00248 \\
		\hline
	\end{tabular}
	\label{tab:Auction_LR_learning_rates}
\end{table}

\section{{Run Times}} \label{appsec:time}

Experiments were implemented on a AMD EPYC 7742 node using a single core, 96 GB. \Cref{tab:time_taken} shows the average time taken per {initialization} for each gradient method and optimization problem.

\begin{table}[!pht]
	\caption{Average {run times} for 1000 steps of each gradient method when minimizing each objective function.}
	\vspace{.1 in}
	\centering
	\begin{tabular}{|c|c|c|c|c|c|c|c|}
		\hline
		Stats Prob & CD & TD & SGD & TSGD & Adam & Adamax & TrAdamax \\
		\hline
		Fermat--Weber & 0.480 & 0.520 & 0.435 & 0.499 & 0.391 & 0.460 & 0.464 \\
		Fr\'echet Mean & 0.347 & 0.363 & 0.301 & 0.342 & 0.294 & 0.422 & 0.370 \\
		\makecell{2-Wasserstein \\ Projection} & 2.739 & 2.134 & 2.392 & 2.209 & 2.341 & 2.246 & 2.272 \\
		\makecell{$\infty$-Wasserstein \\ Projection} & 2.947 & 3.049 & 2.772 & 2.565 & 2.501 & 3.162 & 2.807 \\
		Polytope PCA & 1.117 & 1.081 & 1.014 & 1.052 & 1.037 & 1.174 & 1.029 \\
		Logistic Regression & 0.987 & 1.015 & 0.962 & 0.934 & 0.943 & 1.084 & 0.943 \\
		Linear Regression & 0.494 & 0.499 & 0.386 & 0.420 & 0.386 & 0.537 & 0.433 \\
		\hline
	\end{tabular}
	\label{tab:time_taken}
\end{table}

\pagebreak

\section{Further Figures}\label{appsec:further_figures}

\subsection{Fr\'echet Means}\label{appsubsec:FMs}

\begin{figure}[!pht]
	\centering
	\includegraphics[width = 0.85\textwidth]{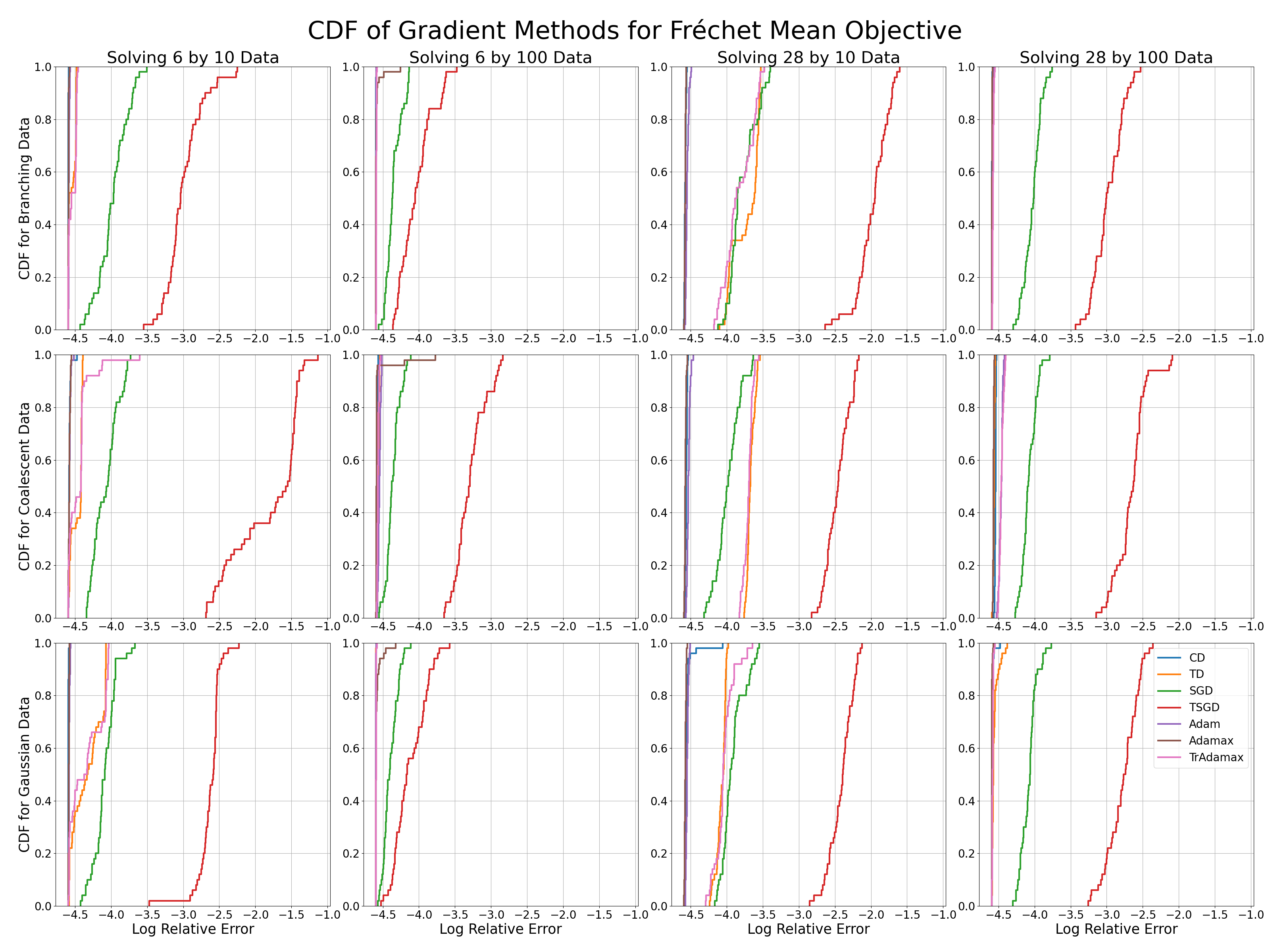}
	\caption{The CDF of the log relative error after 1000 steps for each gradient method with tuned learning rates across 50 random (Gaussian) initializations when minimizing the Fr\'echet mean objective. Each subfigure corresponds to a dataset with sample size $K=10,100$ of dimensionality $N=6,28$, sampled from a branching process, coalescent process, or a Gaussian distribution.}
	\label{fig:FM_CDF_plot}
\end{figure}

Here we include the results for experiments run for the computation of Fr\'echet means. (\cref{tab:FM_mean_log_errors}, \cref{fig:FM_CDF_plot})

\begin{table}[!pht]
	\caption{The mean log relative error of each gradient method after 1000 steps with tuned learning rates for branching process, coalescent process, and Gaussian datasets of size $K=10,100$ in $\TPT{6}$, $\TPT{28}$ when minimizing the Fr\'echet mean objective function. The minimal mean errors for each dataset are in bold.}
	\vspace{.1 in}
	\centering
	\begin{tabular}{|c|c|c|c|c|c|c|c|}
		\hline
		Data & CD & TD & SGD & TSGD & Adam & Adamax & TrAdamax \\
		\hline
		6$\times$10 Branching Data & \textbf{-4.59} & -4.54 & -3.98 & -2.99 & -4.58 & \textbf{-4.59} & -4.53 \\
		6$\times$10 Coalescent Data & \textbf{-4.58} & -4.47 & -4.07 & -1.82 & -4.57 & \textbf{-4.58} & -4.45 \\
		6$\times$10 Gaussian Data & \textbf{-4.59} & -4.33 & -4.09 & -2.63 & -4.58 & \textbf{-4.59} & -4.35 \\
		6$\times$100 Branching Data & \textbf{-4.59} & \textbf{-4.59} & -4.35 & -4.03 & \textbf{-4.59} & -4.58 & \textbf{-4.59} \\
		6$\times$100 Coalescent Data & \textbf{-4.57} & -4.56 & -4.36 & -3.29 & -4.55 & -4.56 & \textbf{-4.57} \\
		6$\times$100 Gaussian Data & -4.59 & \textbf{-4.60} & -4.39 & -4.12 & -4.59 & -4.58 & -4.59 \\
		28$\times$10 Branching Data & \textbf{-4.58} & -3.73 & -3.78 & -1.98 & -4.55 & \textbf{-4.58} & -3.83 \\
		28$\times$10 Coalescent Data & -4.56 & -3.67 & -3.98 & -2.46 & -4.53 & \textbf{-4.57} & -3.70 \\
		28$\times$10 Gaussian Data & -4.55 & -4.07 & -3.91 & -2.41 & -4.55 & \textbf{-4.57} & -4.04 \\
		28$\times$100 Branching Data & \textbf{-4.59} & -4.58 & -4.03 & -2.98 & -4.58 & \textbf{-4.59} & -4.57 \\
		28$\times$100 Coalescent Data & -4.54 & -4.56 & -4.08 & -2.66 & -4.47 & \textbf{-4.57} & -4.46 \\
		28$\times$100 Gaussian Data & \textbf{-4.59} & -4.56 & -4.08 & -2.78 & \textbf{-4.59} & \textbf{-4.59} & -4.58 \\
		\hline
	\end{tabular}
	\label{tab:FM_mean_log_errors}
\end{table}

\subsection{Wasserstein Projections}\label{appsubsec:WDs}

From \Cref{tab:WD2_mean_log_errors,tab:WDinf_mean_log_errors}, we see that when solving the 2-Wasserstein projection problem, tropical methods are optimal in all cases other than $28\times10$ coalescent data. For the $\infty$-Wasserstein projection problem in small dimensions, tropical descent outperforms other methods by some margin, while the comparative performances are less consistent for the $\infty$-Wasserstein projection problem in high dimensions.

\begin{table}[!pht]
	\caption{The mean log relative error of each gradient method after 1000 steps with tuned learning rates for branching process, coalescent process, and Gaussian datasets of size $K=10,100$ in $\TPT{6}$, $\TPT{28}$ when minimizing the $2$-Wasserstein projection objective function. The minimal mean errors for each dataset are in bold.}
	\vspace{.1 in}
	\centering
	\begin{tabular}{|c|c|c|c|c|c|c|c|}
		\hline
		Data & CD & TD & SGD & TSGD & Adam & Adamax & TrAdamax \\
		\hline
		6$\times$10 Branching Data & -3.02 & \textbf{-4.09} & -2.16 & -2.11 & -2.95 & -3.00 & -3.89 \\
		6$\times$10 Coalescent Data & -3.71 & -3.92 & -2.44 & -3.13 & -3.78 & -3.68 & \textbf{-4.12} \\
		6$\times$10 Gaussian Data & -3.59 & -3.83 & -2.87 & -2.75 & -3.61 & -3.57 & \textbf{-3.89} \\
		6$\times$100 Branching Data & -3.44 & \textbf{-4.24} & -2.91 & -1.60 & -3.79 & -3.79 & -4.02 \\
		6$\times$100 Coalescent Data & -4.26 & \textbf{-4.55} & -3.27 & -2.58 & -4.39 & -4.45 & -4.50 \\
		6$\times$100 Gaussian Data & -3.69 & -4.17 & -3.33 & -2.73 & -3.78 & -3.78 & \textbf{-4.18} \\
		28$\times$10 Branching Data & -3.39 & -4.31 & -3.16 & -2.03 & -3.53 & -3.55 & \textbf{-4.32} \\
		28$\times$10 Coalescent Data & -4.11 & -3.79 & -3.39 & -2.42 & \textbf{-4.17} & -4.12 & -3.80 \\
		28$\times$10 Gaussian Data & -3.43 & \textbf{-3.91} & -3.21 & -2.15 & -3.50 & -3.44 & -3.83 \\
		28$\times$100 Branching Data & -3.61 & \textbf{-4.55} & -3.37 & -2.65 & -3.85 & -3.85 & -4.51 \\
		28$\times$100 Coalescent Data & -4.11 & \textbf{-4.57} & -3.89 & -3.32 & -4.21 & -4.19 & -4.56 \\
		28$\times$100 Gaussian Data & -3.66 & \textbf{-4.52} & -3.41 & -2.61 & -3.82 & -3.80 & -4.49 \\
		\hline
	\end{tabular}
	\label{tab:WD2_mean_log_errors}
\end{table}

\begin{table}[!pht]
	\caption{The mean log relative error of each gradient method after 1000 steps with tuned learning rates for branching process, coalescent process, and Gaussian datasets of size $K=10,100$ in $\TPT{6}$, $\TPT{28}$ when minimizing the $\infty$-Wasserstein projection objective function. The minimal mean errors for each dataset are in bold.}
	\vspace{.1 in}
	\centering
	\begin{tabular}{|c|c|c|c|c|c|c|c|}
		\hline
		Data & CD & TD & SGD & TSGD & Adam & Adamax & TrAdamax \\
		\hline
		6$\times$10 Branching Data & -1.63 & \textbf{-4.32} & -0.57 & -0.91 & -1.68 & -1.69 & \textbf{-4.32} \\
		6$\times$10 Coalescent Data & -4.00 & \textbf{-4.47} & -0.99 & -1.43 & -4.03 & -4.04 & -4.46 \\
		6$\times$10 Gaussian Data & -3.77 & \textbf{-4.42} & -1.76 & -3.33 & -3.76 & -3.75 & -4.40 \\
		6$\times$100 Branching Data & -3.35 & \textbf{-4.48} & -1.94 & -0.95 & -3.34 & -3.36 & -4.47 \\
		6$\times$100 Coalescent Data & -3.85 & \textbf{-4.51} & -2.06 & -2.21 & -3.84 & -3.85 & \textbf{-4.51} \\
		6$\times$100 Gaussian Data & -3.63 & \textbf{-4.40} & -1.61 & -1.22 & -3.70 & -3.63 & -4.36 \\
		28$\times$10 Branching Data & -3.68 & -3.53 & -2.22 & -1.17 & \textbf{-4.11} & -3.93 & -3.17 \\
		28$\times$10 Coalescent Data & -4.22 & -3.95 & -2.76 & -1.92 & \textbf{-4.42} & -4.36 & -3.94 \\
		28$\times$10 Gaussian Data & -4.02 & -3.14 & -2.13 & -0.95 & \textbf{-4.43} & -4.30 & -2.85 \\
		28$\times$100 Branching Data & -3.46 & -3.58 & -1.37 & -1.19 & \textbf{-3.68} & -3.59 & -3.39 \\
		28$\times$100 Coalescent Data & -4.41 & -3.92 & -2.09 & -2.00 & \textbf{-4.54} & -4.51 & -3.81 \\
		28$\times$100 Gaussian Data & -2.98 & \textbf{-3.30} & -1.64 & -1.61 & -3.07 & -3.04 & -2.80 \\
		\hline
	\end{tabular}
	\label{tab:WDinf_mean_log_errors}
\end{table}

\clearpage

\bibliographystyle{authordate3}

\bibliography{Sources.bib}

\end{document}